\documentclass[11pt,a4paper]{amsart} 
\usepackage{amsmath,amssymb,amsthm,enumerate,graphicx} 
\usepackage{graphicx} 
\usepackage[all]{xy} 
\usepackage[hypertex]{hyperref} 
 
\numberwithin{equation}{section} 
 
\theoremstyle{plain} 
\newtheorem{theorem} {Theorem} [section] 
\newtheorem*{thm} {Theorem} 
\newtheorem{lemma} [theorem] {Lemma} 
\newtheorem{corollary} [theorem] {Corollary} 
\newtheorem{proposition} [theorem] {Proposition}

\theoremstyle{definition} 
\newtheorem{definition} [theorem] {Definition} 
\newtheorem{remark} [theorem] {Remark}

\renewcommand \parallel {/\kern-3pt/}

\newcommand \N {\mathbb N} 
\newcommand \Z {\mathbb Z}

\newcommand \Esp {\rule[-7pt]{0pt}{10pt}}

\newcommand \Hom {\operatorname{Hom}}

\newcommand \kk {\mathrm{k}} 

\title[Cohomology of quiver algebras]{Comparison morphisms and the \\ 
Hochschild cohomology ring of \\ 
truncated quiver algebras} 
\author{Guillermo Ames} 
\address{Dpto. Matem\'atica, Facultad de Ciencias Exactas, UBA} 
\email{lgames@dm.uba.ar} 
 
\author{Leandro Cagliero} 
\address{CIEM-FaMAF, Universidad Nacional de C\'ordoba} 
\email{cagliero@famaf.unc.edu.ar} 
 
\author{Paulo Tirao} 
\address{CIEM-FaMAF, Universidad Nacional de C\'ordoba} 
\email{ptirao@famaf.unc.edu.ar} 
\thanks{Partially supported by CONICET and Secyt-UNC from Argentina} 

\begin{document} 
 
\begin{abstract} 
 
A main contribution of this paper is the explicit construction of 
comparison morphisms between the standard bar resolution and 
Bardzell's minimal resolution for truncated quiver algebras (TQA's). 
 
As a direct application we describe explicitely the Yoneda product and derive 
several results on the structure of the cohomology ring of TQA's. 
For instance, we show that the product of odd degree cohomology classes is always zero. 
We prove that TQA's associated with quivers with no cycles or with neither sinks nor sources 
have trivial cohomology rings. 
On the other side we exhibit a fundamental example of a TQA with non trivial cohomology ring. 
Finaly, for truncated polyniomial algebras in one variable, we construct explicit cohomology classes 
in the bar resolution and give a full description of their cohomology ring. 
 
\end{abstract} 
\maketitle

\section{Introduction} 
\label{sec:intro} 
 
To any finite quiver $\Delta$ and any field $\kk$ one associates a 
$\kk$-algebra ${\mathrm{k}}\Delta$, 
the {\sl path algebra} or {\sl quiver algebra} of $\Delta$, where the set of 
vertices $\Delta_0$ and the sets of $k$-paths $\Delta_k$ form a ${\mathrm{k}}$-basis 
and the product is given by concatenation of paths (see \S \ref{sec:Preliminaries}). 
 
Quiver algebras and their quotients arise 
in many contexts and have been extensively studied. 
A result of Gabriel \cite{G} establishes that for every finite dimensional 
${\mathrm{k}}$-algebra $A$ such that 
 $A/r={\mathrm{k}}\times\dots\times{\mathrm{k}}$, 
 where $r$ is the Jacobson radical of $A$, 
there exists a finite quiver $\Delta$, the Gabriel quiver of $A$, 
and an epimorphism $\varphi:{\mathrm{k}}\Delta\to A$ such that 
$(\Delta_N)\subset\ker\varphi\subset(\Delta_2)$ for some $N\ge2$. 
Here $(\Delta_k)$ is the two-sided ideal generated by $\Delta_k$. 
 
{\sl Monomial algebras} are those for which $\ker\varphi$ is generated by monomials. 
In the particular case when $\ker\varphi=(\Delta_N)$, the algebra $A$ is a 
{\sl truncated quiver algebra}, denoted TQA from now on. 
This property of $A$ turns out to be intrinsic \cite{Ci2} which makes 
TQA's a distinguished class. 
 
For these classes of algebras Bardzell \cite{Ba} introduced a minimal resolution that 
plays a key role in the treatment of homological questions and problems. 
 
A main contribution of this paper is the explicit construction of 
comparison morphisms between the standard bar resolution and 
Bardzell's minimal resolution for TQA's. 
We believe that such morphisms, sought for a long time, should have 
many applications. 
Our construction in this case could inspire others to find comparison 
morphisms for wider classes of algebras, hopefully for all monomial algebras. 
 
As a direct application we describe explicitely the Yoneda product and derive 
several results on the structure of the cohomology ring of TQA's. 
In the near future we will complete a full description of this ring.

\subsection{A brief account of known results} 

Since the early nineties several authors have investigated 
different homological questions for a number of classes of monomial 
algebras, including TQA's. 
 
In the context of truncated quiver algebras, Cibils \cite{Ci1} 
proved that $H_n(A,A) = 0$ for all $n>0$ if the quiver has no 
oriented cycles. A shorter proof of this fact was later given in 
\cite{Ci2}. On the other hand, Liu and Zhang \cite{LZ} showed 
that $H_n(A,A) = 0$ for all $n>0$ if and only if the quiver has no 
oriented cycles of some specific lengths. More recently, 
Sk{\"o}ldberg \cite{Sk1} gave a complete description of the 
homology of TQA's. His computations are based 
on the use of Bardzell's minimal resolution. In the same paper he treats 
also the case of quadratic monomial algebras with an analogous 
approach. 
 
The resolution used in \cite{Sk1} was introduced by Bardzell \cite{Ba} for 
monomial algebras, and has shown to be an efficient tool for computations 
in contrast to the usual bar resolution. 
Recently, N.\ Marconnet in \cite{M} constructed a comparison morphism between the two 
for the first non trivial degree, in the context of cubic Artin-Schelter regular algebras. 
 
The first cohomology computations for TQA's 
appeared in \cite{Ci2} where the second cohomology group is 
described to study formal deformations and to characterize 
rigidity. In the subsequent paper \cite{Ci3} these 
results were extended to the class of monomial algebras. 
 
A description of the whole cohomology of TQA's, over fields of 
characteristic zero, was given by Locateli in \cite{Lo}. 
Her computations also rely on the use of 
Bardzell's minimal resolution, and cohomology classes are 
represented by pairs of parallel paths. The particular case of 
truncated cycle algebras is treated separately. 
Recently in \cite{XHJ} the case of arbitrary characteristic was solved. 
 
The determination of the structure of the full cohomology ring is 
still a difficult problem that has been addressed in a number of 
cases. 
 
For instance, for a radical square zero algebra, a description of 
the Yoneda product on Hochschild cohomology is given in \cite{Ci4} 
and it is shown that this algebra is finitely generated only for 
the case when the underlying quiver is a cycle or it has no oriented 
cycles. 
 
For truncated cycle algebras, the complete structure of the 
cohomology ring was determined in \cite{BLM} and independently in 
\cite{EH}, showing in particular that the Yoneda product is non 
trivial and the cohomology ring is finitely generated. Cycle 
algebras are examples of self-injective Nakayama algebras. In 
\cite{BLM} the authors present in contrast some examples of non 
injective Nakayama algebras for which the product is trivial (in 
non zero degree) and in particular the cohomology ring, which is 
infinite, is not finitely generated. 
 
Being $H^*(A,A)$ a graded commutative $\kk$-algebra every 
homogeneous element of odd degree squares to zero ($char\ \kk\ne 
2$) and if $\mathcal{N}$ is the ideal generated by the homogeneous nilpotent 
elements, then $H^*(A,A)/\mathcal{N}$ is a commutative $\kk$-algebra. One 
expects to gain information for the full cohomology ring from this 
simpler one. 
 
In \cite{SS} it was conjectured that $H^*(A,A)/\mathcal{N}$ is finitely 
generated as $\kk$-algebra for any finite dimensional algebra $A$. 
This was recently proved in \cite{GSS} for monomial algebras and 
was already known for other classes (see \cite{GSS}). In \cite{GS} 
the quotient $H^*(A,A)/\mathcal{N}$ was determined for the subclass of 
stacked monomial algebras, class that contains TQA's. These 
results applied to TQA's yield $H^*(A,A)/\mathcal{N}\simeq \kk$. However, 
since $H^*(A,A)$ is, in general, infinite dimensional over $\kk$, the structure of 
the full ring $H^*(A,A)$ remains open.

\subsection{An overview of the main results} 
 
Given $A$ an associative $\kk$-algebra with unit, 
the Hochschild cohomology groups $H^n(A,A)$ are, by definition, 
the groups $\text{Ext}_{A^e}^n(A,A)$ where $A^e=A\otimes_{\kk} A^{op}$. 
The  natural identification between $A$-bimodules and 
left $A^e$-modules gives the definition of projective $A$-bimodule and 
$A$-bimodule homomorphism. 
 
The standard bar resolution of $A$ is the $A^e$-projective resolution 
\[ 
\cdots\rightarrow A\otimes A^{\otimes n}\otimes A \stackrel{b}\rightarrow 
A\otimes A^{\otimes(n-1)}\otimes A\cdots \stackrel{b}\rightarrow A\otimes A 
\otimes A\stackrel{b}\rightarrow A\otimes A\stackrel{\epsilon}\rightarrow A 
\] 
and since 
$ 
\Hom_{A^e}(A\otimes A^{\otimes n}\otimes A,A) 
\simeq\Hom_{\mathrm{k}}(A^{\otimes n},A) 
$, 
the associated Hochschild complex is 
\[ 
A\stackrel{b}\rightarrow\Hom_{\mathrm{k}}(A,A) 
\stackrel{b}\rightarrow\cdots \stackrel{b}\rightarrow 
\Hom_{\mathrm{k}}(A^{\otimes(n-1)},A)\stackrel{b}\rightarrow 
\Hom_{\mathrm{k}}(A^{\otimes n},A)\stackrel{b}\rightarrow\cdots 
\] 
The cohomology group $H^*(A,A)$ 
has a ring structure given by the Yoneda product which 
coincides with the cup product defined as follows. 
Given two cochains, 
\[ 
f\in\Hom_{\mathrm{k}}(A^{\otimes m},A),\quad 
g\in\Hom_{\mathrm{k}}(A^{\otimes n},A) 
\] 
the cup product of $f$ and $g$ is the cochain 
$ 
f\cup g\in\Hom_{\mathrm{k}}(A^{\otimes (m+n)},A) 
$ 
defined by 
\[ f\cup g\;(\alpha_1\otimes\dots\otimes\alpha_{m+n})= 
f(\alpha_1\otimes\dots\otimes\alpha_{m})\,g(\alpha_{m+1}\otimes\dots\otimes\alpha_{m+n}). 
   \] 
 
For TQA's, the bar resolution can be slightly simplified 
with the $A$-bimodule 
\[ 
{\bf Q}_n= 
A\otimes_{\Delta_0}A_+^{\otimes^n_{\Delta_0}}\otimes_{\Delta_0}A, 
\] 
in place of $A\otimes A^{\otimes n}\otimes A$, 
where $A_+$ is the ideal $\oplus_{n=1}^{N-1}{\mathrm{k}}\Delta_n$. 
For the definition of the differential see \S \ref{subsec:reduced-resolution}. 
 
In contrast to this resolution, there is the following minimal resolution 
$({\bf P},d)$, due to Bardzell \cite{Ba} 
(cf.\ \cite{AG}, \cite{BK} and \cite{Ha}), 
where the $A^e$-projective modules are 
\[ 
{\bf P}_n= 
\begin{cases} 
A\otimes_{\Delta_{0}}{\mathrm{k}}\Delta_{kN}\otimes_{\Delta_{0}}A, &\text{if $n=2k$;} \\ 
A\otimes_{\Delta_{0}}{\mathrm{k}}\Delta_{kN+1}\otimes_{\Delta_{0}}A, &\text{if $n=2k+1$.} 
\end{cases} 
\] 
One has 
\begin{equation*} 
{\bf P}^*_n=\Hom_{A^e}({\bf P}_n,A)\simeq 
\begin{cases} 
\Hom_{\Delta_{0}^e}({\mathrm{k}}\Delta_{kN},A), &\text{if $n=2k$;} \\ 
\Hom_{\Delta_{0}^e}({\mathrm{k}}\Delta_{kN+1},A), &\text{if $n=2k+1$.} 
\end{cases} 
\end{equation*} 
 
The definition of the differential is in \S \ref{subsec:minimal-resolution}. 
For more details see \S \ref{sec:cohomology ring}. 
 
We define the following $A^e$-morphisms between these two resolutions in both directions. 
See \S \ref{sec:comparison-morphisms} for complete details. 
 
First, let ${\bf F}:{\bf P}\to{\bf Q}$ be the $A$-bimodule extension of the map defined 
on $p_0=1\otimes v_1\dots v_{kN}\otimes1\in {\bf P}_{2k}$ and $p_1=1\otimes v_1\dots v_{kN+1}\otimes1\in {\bf P}_{2k+1}$ by 
\begin{multline*} 
{\bf F}_{2k}(p_0)= \\ 
\sum 
1 
[ 
\underbrace{v_1\dots v_{x_1}}_{x_1} 
| 
\underbrace{v_{1+x_1}}_{1} 
| 
\underbrace{\dots v_{1+x_1+x_2}}_{x_2} 
| 
\underbrace{v_{2+x_1+x_2}}_{1} 
|\dots\dots| 
\underbrace{v_{k+\sum x_j}}_{1} 
] 
\underbrace{\dots\dots v_{kN}}_{kN-k-\sum x_j}, 
\end{multline*} 
\begin{multline*} 
{\bf F}_{2k+1}(p_1)= \\ 
\sum 
1 
[\! 
\underbrace{\!\!v_1\!\!}_{1}\! 
| 
\underbrace{v_2\dots v_{1+x_1}}_{x_1} 
| 
\underbrace{v_{2+x_1}}_{1} 
| 
\underbrace{\dots v_{2+x_1+x_2}}_{x_2} 
|\dots\dots| 
\underbrace{v_{k\!+\!1+\sum x_j}}_{1} 
] 
\underbrace{\dots\dots v_{kN+1}}_{kN-k-1-\sum x_j}, 
\end{multline*} 
where the sum is taken over all $k$-tuples $(x_1,\dots,x_k)\in\Z^k$ 
such that $1\le x_i<N$.

Second, let ${\bf G}:{\bf Q}\to{\bf P}$ be the $A$-bimodule extension of the map defined, 
on $q=1[\alpha_1|\dots|\alpha_n]1=a_1^1\dots a_{|\alpha_1|}^1a_1^2\dots a_{|\alpha_2|}^2\dots\dots a_1^{n} 
\dots a_{|\alpha_{n}|}^{n}=v_1\dots v_{|q|}\in{\bf Q}_n$, by 
 
\begin{align*} 
{\bf G}_{2k}(q)= 
\begin{cases} 
1\otimes v_1\dots v_{kN}\otimes v_{kN+1}\dots v_{|q|},& 
\text{if }\alpha_{2i-1}\alpha_{2i}=0\text{ for }i=1\dots k; \\ 
0, &\text{otherwise}; 
\end{cases} 
\end{align*} 
 
\begin{align*} 
{\bf G}_{2k+1}(q)= 
\begin{cases} 
\displaystyle\sum_{j=1}^{|\alpha_1|} 
v_1\dots v_{j-1}\otimes v_{j}\dots v_{kN+j}\otimes v_{kN+j+1}\dots v_{|q|},& 
\text{if }\alpha_{2i}\alpha_{2i+1}=0\\[-5mm] 
& \text{for }i=1\dots k; \\[2mm] 
0, &\text{otherwise}. 
\end{cases} 
\end{align*} 
 
Our first result is the following theorem (see Theorem \ref{MainComparison}). 
 
\begin{thm} 
The morphisms ${\bf F}$ and ${\bf G}$ between the resolutions ${\bf P}$ and ${\bf Q}$ 
are both comparison morphisms. 
\end{thm} 
 
The proofs are long, sections  \S \ref{sec:proof-F} and \S \ref{sec:proof-G} are exclusively devoted 
to them. They are subtle and give an insight on the non trivial combinatorics 
underlying this problem. 
 
\bigskip 
 
With the comparison morphisms at hand we describe first the Yoneda 
product at the level of the minimal resolution and then determine the product in cohomology. 
The minimal resolution is naturally bigraded, but the product in this resolution 
is not compatible with this bigrading. 
However, the product at the cohomology level, which is essentially given 
by concatenation of paths, is compatible with the bigrading making 
$H^*(A,A)$ a bigraded ring. 
 
More precisely, let $\vee$ be the product in ${\bf P}^*$ defined in the following way. 
For $(\alpha,\pi)\in {\bf P}^*_{n_1}$ and $(\beta,\tau)\in{\bf P}^*_{n_2}$, 
 
\[(\alpha,\pi)\vee(\beta,\tau)=\begin{cases} 
(\alpha\beta,\pi\tau), & \text{ if $n_1$ or $n_2$ is even} \\ 
0, & \text{ otherwise.}\end{cases}\] 
 
We then have the following result. 
 
\begin{thm} 
Let $A$ be an $N$-TQA. 
Then the $\vee$ product in ${\bf P}^*$ induces the Yoneda product in $H^n(A,A)$ and in particular: 
\begin{enumerate}[(i)] 
\item The product of two odd degree cohomology classes is zero. 
\item If $f_1,...,f_N$ are cohomology classes of positive degree, then $f_1... f_N=0$. 
\item $H^*(A,A)/\mathcal{N}=\kk$, where $\mathcal{N}$ is the ideal generated by homogeneous nilpotent elements. 
\end{enumerate} 
\end{thm} 
 
This result extends the analogous result in \cite{BLM} for truncated cycle algebras. 
Part $(iii)$ can be deduced from the results in \cite{GS}. 
 
This theorem allow us to derive a number of results on the structure of 
the full cohomology ring of TQA's. 
In this paper we investigate under which conditions is the cohomology ring trivial, 
meaning that the subalgebra $\oplus_{n\ge 1} H^n(A,A)$ has trivial product. 
It was believed that, generically, this was the case. 
Nevertheless examples of algebras with non trivial product in cohomology 
appeared in \cite{GMS}. Recently in \cite{GS,GSS} examples within the class of 
monomial algebras are given. 
On the other hand, Bustamante and Gatica \cite{BG} proved that the product is zero for monomial 
algebras with no oriented cycles. 
 
For the class of TQA's we prove in \S\ref{sec:examples} the following theorem. 
 
\begin{thm} 
Let $\Delta$ be a quiver satisfying one of the following conditions. 
 \begin{enumerate}[(i)] 
  \item $\Delta$ has no oriented cycles. 
  \item $\Delta$ is not an oriented cycle and has neither sinks nor sources. 
 \end{enumerate} 
Then the subalgebra $\oplus_{n\ge 1}H^n(A,A)$ with the Yoneda product is 
trivial. 
\end{thm} 
 
On the other direction, we consider 
the cohomology ring of TQA'a associated with the quiver 
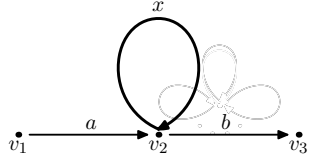
\begin{center} 
   \includegraphics[width=40mm]{quiver.2} 
\end{center} 
and prove in \S\ref{subsec:example} the following result. 
\begin{thm} 
Let $A$ be an $N$-TQA associated with the above quiver $\Delta$. 
Then, for all $n\in\mathbb{N}$, there exist non zero cohomology classes $\omega_{n,j}\in H^n(A,A)$, 
$j=1,\dots, N-1$, such that 
\[ 
\omega_{n_1,j_1}\cup\omega_{n_2,j_2}= 
\begin{cases} 
\omega_{n_1+n_2,j_1+j_2}, &\text{ if $n_1$ or $n_2$ is even and $j_1+j_2<N$}; \\ 
0,&\text{  otherwise}. 
\end{cases} 
\] 
\end{thm}

This theorem gives many examples of TQA's containing loops (and thus oriented cycles) whose cohomology 
ring contain nilpotent elements that are factorized as a product of two other nilpotent elements. 
On a full description of the cohomology ring of arbitrary TQA's this example should play 
a fundamental role. 
 
At the end of the paper we use the comparison morphisms to construct explicit cohomology classes 
in the bar resolution. We give a full description of the cohomology ring of 
truncated polynomial algebras in one variable. 
 
{\bf Acknowledgements.} 
We wish to thank heartily Marco Farinati for very helpful conversations on 
various matters in homological algebra. 
We also thank Victor Ginzburg for his interest 
in this work and for calling to our attention the work \cite{M}.

\section{Preliminaries} 
\label{sec:Preliminaries} 

\subsection{Quiver algebras} 
 
Let $\Delta$ be a finite quiver, that is a finite directed graph in which 
multiple arrows and loops are allowed. 
In this paper all quivers shall be assumed to be finite and connected. 
 
The set of vertices and arrows of $\Delta$ are denoted by 
$\Delta_0$ and $\Delta_1$ respectively. 
To each arrow $a\in\Delta_1$ we 
associate its source vertex $o(a)$, and its end vertex $t(a)$. 
A path $\alpha$ is a sequence of arrows 
$\alpha=a_1\dots a_n$ such that $t(a_i)=o(a_{i+1})$. 
The length $|\alpha|$ of a path $\alpha$ is the number of 
arrows of it and the set of paths of length $n$ is denoted by $\Delta_n$. 
We find it convenient to consider the vertices as paths of length zero. 
For a path $\alpha=a_1\dots a_n\in\Delta_n$, 
we set $o(\alpha)=o(a_{1})$ and $t(\alpha)=t(a_n)$. 
 
Let ${\mathrm{k}}$ be any field of characteristic 0. 
Let ${\mathrm{k}}\Delta_n$ be the ${\mathrm{k}}$-vector space with basis 
$\Delta_n$ and let 
${\mathrm{k}}\Delta=\bigoplus_{n\ge0}{\mathrm{k}}\Delta_n$. 
The quiver algebra associated to $\Delta$ is ${\mathrm{k}}\Delta$ 
with multiplication given by concatenation of paths. 
If $\alpha=a_1\dots a_m\in\Delta_m$ and $\beta=b_1\dots b_n\in\Delta_n$, then 
$\alpha\beta=a_1\dots a_mb_1\dots b_n\in\Delta_{m+n}$, if $t(\alpha)=o(\beta)$, 
or zero otherwise. 
It is clear that ${\mathrm{k}}\Delta$ is a graded algebra with 
unit $1=\sum_{p\in\Delta_0}p$ and degree $n$ component ${\mathrm{k}}\Delta_n$. 
 
A \textsl{truncated quiver algebra} $A$ is a quotient 
$A={\mathrm{k}}\Delta/I^N$, where $I$ is the ideal generated by $\Delta_1$ 
and $N$ is a positive integer. 
Since $I^N$ is an homogeneous ideal, truncated quiver algebras are graded. 
 
Given a truncated quiver algebra $A$, 
we shall make no distinction between an element 
$\alpha\in\bigoplus_{n=0}^{N-1}{\mathrm{k}}\Delta_n\subset {\mathrm{k}}\Delta$ 
and its quotient projection in $A$. 
In particular, the set 
\[ 
\mathcal B=\bigcup_{n=0}^{N-1}\Delta_n 
\] 
 is a ${\mathrm{k}}$-basis of $A$. 
 
We finally point out that elements $\alpha,\beta\in A$ 
will frequently be assumed to be in $\mathcal B$ and, in these cases, 
$\alpha=a_1\dots a_{|\alpha|}$ or $\beta=b_1\dots b_{|\beta|}$ 
will be their arrow decomposition. 
 
\subsection{Hochschild cohomology} 
Given an associative ${\mathrm{k}}$-algebra with unit $A$, 
the Hochschild cohomology groups $H^*(A,A)$ of $A$ with coefficients in the 
$A$-bimodule $A$ are, by definition, 
$\text{Ext}_{A^e}^*(A,A)$ where $A^e=A\otimes_{\kk} A^{op}$. 
The  natural identification between $A$-bimodules and 
left $A^e$-modules gives the definition of projective $A$-bimodule and 
$A$-bimodule homomorphism. 
 
We recall that the standard bar resolution of a ${\mathrm{k}}$-algebra with unit $A$ 
is the $A^e$-projective resolution of $A$ 
\[ 
\cdots\rightarrow A\otimes A^{\otimes n}\otimes A \stackrel{b}\rightarrow A\otimes A^{\otimes(n-1)} 
\otimes A\cdots \stackrel{b}\rightarrow A\otimes A 
\otimes A\stackrel{b}\rightarrow A\otimes A\stackrel{\epsilon}\rightarrow A 
\] 
where $\epsilon(\alpha\otimes\beta)=\alpha\beta$ and the differential $b$ in degree $n$ 
is given by 
\begin{multline*} 
b_n(\alpha_0\otimes \alpha_1\dots \alpha_n\otimes \alpha_{n+1})= 
\alpha_{0}\alpha_1\otimes \alpha_2\dots \alpha_n\otimes \alpha_{n+1}\ + \\ 
\sum_{i=1}^{n-1}(-1)^i \alpha_{0}\otimes \alpha_1\dots(\alpha_i\alpha_{i+1})\dots \alpha_n\otimes \alpha_{n+1} + 
(-1)^n \alpha_{0}\otimes \alpha_1\dots \alpha_{n-1}\otimes \alpha_n \alpha_{n+1}. 
\end{multline*} 
Since 
$ 
\Hom_{A^e}(A\otimes A^{\otimes n}\otimes A,A) 
\simeq\Hom_{\mathrm{k}}(A^{\otimes n},A) 
$, 
the associated Hochschild complex is 
\[ 
A\stackrel{b}\rightarrow\Hom_{\mathrm{k}}(A,A) 
\stackrel{b}\rightarrow\cdots \stackrel{b}\rightarrow 
\Hom_{\mathrm{k}}(A^{\otimes(n-1)},A)\stackrel{b}\rightarrow 
\Hom_{\mathrm{k}}(A^{\otimes n},A)\stackrel{b}\rightarrow\cdots 
\] 
The cohomology group $H^*(A,A)$ 
has a ring structure given by the Yoneda product which 
coincides with the cup product defined as follows. 
The cup product is graded commutative. 
Given two cochains, 
\[ 
f\in\Hom_{\mathrm{k}}(A^{\otimes m},A),\quad 
g\in\Hom_{\mathrm{k}}(A^{\otimes n},A) 
\] 
the cup product of $f$ and $g$ is the cochain 
$ 
f\cup g\in\Hom_{\mathrm{k}}(A^{\otimes (m+n)},A) 
$ 
defined by 
\[ f\cup g\;(\alpha_1\otimes\dots\otimes\alpha_{m+n})= 
f(\alpha_1\otimes\dots\otimes\alpha_{m})\,g(\alpha_{m+1}\otimes\dots\otimes\alpha_{m+n}). 
   \] 
We finally recall that the Hochschild cohomology of the direct sum of 
two $\mathrm{k}$-algebras is the direct sum of their Hochschild cohomologies. 
Thus we shall restrict ourselves to finite connected quivers.

\section{Two projective resolutions} 
\label{sec:projective-resolutions} 
 
\subsection{The (reduced) bar resolution $({\bf Q},b)$ }\label{Res Q} 
\label{subsec:reduced-resolution} 
 
When $A$ is a truncated quiver algebra the bar resolution given above 
can be slightly simplified by tensoring over ${\mathrm{k}}\Delta_0$, 
as done in \cite{Ci2}. 
More precisely, let 
us denote by $A_+$ the ideal $\oplus_{n=1}^{N-1}{\mathrm{k}}\Delta_n$ 
of $A$ and let 
\[ 
{\bf Q}_n= 
A\otimes_{\Delta_0}A_+^{\otimes^n_{\Delta_0}}\otimes_{\Delta_0}A, 
\] 
$\epsilon(\alpha\otimes\beta)=\alpha\beta$ and for $n>0$ let 
\begin{multline*} 
b_n(\alpha_0[\alpha_1|\dots|\alpha_n]\alpha_{n+1})= 
\alpha_0\alpha_1[\alpha_2|\dots|\alpha_n]\alpha_{n+1}+ \\ 
 + \sum_{i=1}^{n-1}(-1)^i \alpha_0[\alpha_1|\dots|\alpha_i\alpha_{i+1}|\dots| \alpha_n]\alpha_{n+1} + 
(-1)^n \alpha_0[\alpha_1|\dots|\alpha_{n-1}]\alpha_n\alpha_{n+1}. 
\end{multline*} 
Here we use the bar notation 
\[ 
\alpha_0[\alpha_1|\dots|\alpha_n]\alpha_{n+1}= 
\alpha_0\otimes_{\Delta_0} \alpha_1\otimes_{\Delta_0}\dots 
\otimes_{\Delta_0} \alpha_n\otimes_{\Delta_0} \alpha_{n+1}\in{\bf Q}_n. 
\] 
It is not difficult to see that ${\bf Q}_n$ is $A^e$-projective, 
that $b$ is well defined and $b^2=0$ (see \cite{Ci2}). 
 
A ${\mathrm{k}}$-basis of ${\bf Q}_n$ is 
\[ 
\mathcal B_{{\bf Q}_n}=\left\{\alpha_0[\alpha_1|\dots|\alpha_n]\alpha_{n+1}\left| 
\begin{array}{ll} 
i) &\!\!\alpha_j\in\mathcal B\text{ for all $j$; } 
    |\alpha_j|\ge1,\text{ for }j\!=\!1\!\dots\! n\rule[-7pt]{0pt}{10pt} \\ 
ii) &\!\!t(\alpha_j)=o(\alpha_{j+1}),\text{ for } j=0\dots n 
\end{array}\right.\! 
\right\}. 
\] 
Let 
\[ 
\mathcal B'_{{\bf Q}_n}=\left\{1[\alpha_1|\dots|\alpha_n]1\left| 
\begin{array}{ll} 
i) &\alpha_j\in\mathcal B\text{ and } 
    |\alpha_j|\ge1,\text{ for }j=1\dots n\rule[-7pt]{0pt}{10pt} \\ 
ii) &t(\alpha_j)=o(\alpha_{j+1}),\text{ for } j=1\dots n-1 
\end{array}\right.\! 
\right\}. 
\] 
Since 
\[ 
1[\alpha_1|\dots|\alpha_n]1=o(\alpha_{1})[\alpha_1|\dots|\alpha_n]t(\alpha_n) 
\] 
for every element of $\mathcal B'_{{\bf Q}_n}$, it follows that 
$\mathcal B'_{{\bf Q}_n}\subset\mathcal B_{{\bf Q}_n}$ and that 
the set $\mathcal B'_{{\bf Q}_n}$ generates ${\bf Q}_n$ as $A$-bimodule.

As in the case of the bar resolution, 
it is straightforward to check that the map $s:{\bf Q}_n\to{\bf Q}_{n+1}$ defined by 
\[ 
s_n(\alpha_0[\alpha_1|\dots|\alpha_n]\alpha_{n+1})= 
\begin{cases} 
1[\alpha_0|\dots|\alpha_n]\alpha_{n+1}, &\text{if }|\alpha_0|>0; \\ 
0,                       &\text{if }|\alpha_0|=0; 
\end{cases} 
\] 
is a ${\mathrm{k}}$-linear chain contraction of the identity, that is $sb+bs=1$. 
This shows that the complex $({\bf Q},b)$ is exact.

\subsection{The minimal resolution $({\bf P},d)$}\label{Res P} 
\label{subsec:minimal-resolution} 
The Hochschild homology of truncated quiver  algebras $A$ was 
computed by Sk\"oldberg in \cite{Sk1} and the Hochschild cohomology 
was computed by Locateli in \cite{Lo}. 
In both papers, the authors used the minimal $A^e$-projective resolution 
${\bf P}$ of $A$ that we describe below. 
This minimal resolution was introduced in several earlier papers 
(see for instance \cite{AG}, \cite{Ba}, \cite{BK} and \cite{Ha}). 
 
Let 
\[ 
{\bf P}_n= 
\begin{cases} 
A\otimes_{\Delta_{0}}{\mathrm{k}}\Delta_{kN}\otimes_{\Delta_{0}}A, &\text{if $n=2k$;} \\ 
A\otimes_{\Delta_{0}}{\mathrm{k}}\Delta_{kN+1}\otimes_{\Delta_{0}}A, &\text{if $n=2k+1$.} 
\end{cases} 
\] 
In order to simplify the notation, the symbol $\otimes$ 
will always mean $\otimes_{\Delta_{0}}$ for 
elements in ${\bf P}$. 
Let $\epsilon(\alpha\otimes\beta)=\alpha\beta$ and, for $n>0$, let 
$d_n:{\bf P}_n\to{\bf P}_{n-1}$ be defined by 
\begin{align*} 
d_{2k}(\alpha \!\otimes\! v_1\dots v_{kN}\!\otimes\! \beta)&= 
\sum_{j=0}^{N-1} 
\alpha v_1\dots v_j\otimes 
\underbrace{v_{j+1}\dots v_{t}}_{(k-1)N+1}\otimes v_{t+1}\dots v_{kN}\beta \\ 
&= 
\alpha \otimes v_{1}\dots v_{(k-1)N+1}\otimes v_{(k-1)N+2}\dots v_{kN} \beta \Esp\\ 
&\hspace{2cm}+\dots + 
\alpha  v_1\dots v_{N-1}\otimes v_{N}\dots v_{kN}\!\otimes\! \beta, \Esp\\ 
d_{2k+1}(\alpha \!\otimes\! v_1\dots v_{kN+1}\!\otimes\! \beta)&= 
\alpha v_1\!\otimes\! v_{2}\dots v_{kN+1}\!\otimes\! \beta\,-\, 
\alpha \!\otimes\! v_{1}\dots v_{kN}\!\otimes\! v_{kN\!+\!1} \beta. 
\end{align*} 
In particular 
\begin{align*} 
d_{1}(\alpha\otimes v\otimes \beta)&=\alpha v\otimes \beta\,-\,\alpha\otimes v \beta, \\ 
d_{2}(\alpha\otimes v_1\dots v_{N}\otimes \beta)&= 
\sum_{j=0}^{N-1} 
\alpha v_1\dots v_j\otimes v_{j+1}\otimes v_{j+2}\dots v_{N}\beta. 
\end{align*} 
 
Again, it is easy to see that ${\bf P}_n$ is $A^{e}$-projective, 
$d$ is well defined, $d^2=0$ and the set 
\[ 
\mathcal B'_{{\bf P}_n}\!=\left\{1\!\otimes\! v_1\dots v_s\!\otimes\!1\left| 
\begin{array}{ll} 
i)  &\!\!\! s\!=\!kN,\text{ if }n\!=\!2k; 
                  \text{ or }s\!=\!kN\!+\!1,\text{ if }n\!=\!2k\!+\!1 \\ 
ii) &\!\!\! v_i\in\Delta_1\text{ for }i=1\dots s 
\end{array}\right.\!\! 
\right\}. 
\] 
generates ${\bf P}_n$ as $A$-bimodule. 
 
The exactness of $({\bf P},d)$ is not obvious and a 
proof using a spectral sequence argument can be found in \cite{Sk1}. 
Alternatively, we give a chain contraction of the identity 
$r_n:{\bf P}_n\to{\bf P}_{n+1}$ 
in the following proposition. 
 
\begin{proposition} 
Let $r_n:{\bf P}_n\to{\bf P}_{n+1}$ be the $\mathrm{k}$-linear 
map defined on basis elements as follows 
\begin{align*} 
r_{2k}(\alpha\otimes v_1\dots v_{kN}\otimes\beta)&= 
\sum_{j=1}^{|\alpha|} 
a_1\dots a_{j\!-\!1}\otimes 
\underbrace{a_{j}\dots a_{|\alpha|}v_1\dots v_{t}}_{kN+1} 
\otimes v_{t+1}\dots v_{kN}\beta \\ 
&= 
1\otimes \alpha v_1\dots v_{kN-|\alpha|+1}\otimes v_{kN-|\alpha|+2}\dots v_{kN}\beta\Esp \\ 
&\hspace{15mm}+\dots +a_1\dots a_{|\alpha|-1}\otimes a_{|\alpha|}v_1\dots v_{kN}\otimes\beta;\Esp \\ 
r_{2k\!+\!1}(\alpha \otimes\! v_1\dots v_{kN\!+\!1}\!\otimes\beta)&= 
\begin{cases} 
1\otimes \alpha v_{1}\dots v_{kN+1}\otimes\beta,&\text{ if $|\alpha |=N-1$;} \\ 
0, &\text{ if $|\alpha |<N-1$.} 
\end{cases} 
\end{align*} 
Then $rd+dr=1$ and therefore $({\bf P},d)$ is exact. 
\end{proposition} 
\begin{remark} 
For $n=0$ and $n=1$ we have 
\begin{align*} 
r_{0}(\alpha \otimes\beta)&=\sum_{j=1}^{|\alpha |}a_1\dots a_{j-1}\otimes 
a_{j}\otimes a_{j+1}\dots a_{|\alpha|}\beta \\ 
&=1\otimes a_1\otimes a_{2}\dots a_{|\alpha|}\beta 
+\dots +a_1\dots a_{|\alpha|-1}\otimes a_{|\alpha |}\otimes \beta\Esp \\ 
r_{1}(\alpha \otimes v\otimes \beta)&= 
\begin{cases} 
1\otimes \alpha v\otimes \beta,&\text{ if $|\alpha |=N-1$;} \\ 
0, &\text{ if $|\alpha |<N-1$.} 
\end{cases} 
\end{align*} 
\end{remark} 
 
\begin{proof} 
For $n=2k$ we have 
\begin{align*} 
d_{2k+1}&r_{2k}(\alpha\otimes v_1\dots v_{kN}\otimes \beta)\\ 
&= 
\sum_{j=1}^{|\alpha|} 
d_{2k+1}(a_1\dots a_{j-1}\otimes 
\underbrace{a_{j}\dots a_{|\alpha|}v_1\dots v_{t}}_{kN+1} 
\otimes v_{t+1}\dots v_{kN}\beta) \\ 
&= 
\sum_{j=1}^{|\alpha|} 
a_1\dots a_{j}\otimes 
\underbrace{a_{j+1}\dots a_{|\alpha|}v_1\dots v_{t}}_{kN} 
\otimes v_{t+1}\dots v_{kN}\beta \\ 
&\hspace{2cm}- 
\sum_{j=1}^{|\alpha|} 
a_1\dots a_{j-1}\otimes 
\underbrace{a_{j}\dots a_{|\alpha|}v_1\dots v_{t-1}}_{kN} 
\otimes v_{t}\dots v_{kN}\beta \\ 
&= 
\alpha\otimes v_1\dots v_{kN}\otimes \beta \;-\; 
1\otimes \alpha v_1\dots v_{kN-|\alpha|} 
\otimes v_{kN-|\alpha|+1}\dots v_{kN}\beta; 
\end{align*} 
and 
\begin{align*} 
r_{2k-1}&d_{2k}(\alpha\otimes v_1\dots v_{kN}\otimes \beta)\\ 
&= 
\sum_{j=0}^{N-1} 
r_{2k-1}(\alpha v_1\dots v_{j}\otimes 
\underbrace{v_{j+1}\dots v_{t}}_{(k-1)N+1} 
\otimes v_{t+1}\dots v_{kN}\beta) \\ 
&= 
r_{2k-1}(\alpha v_1\dots v_{N-1-|\alpha|}\otimes 
v_{N-|\alpha|}\dots v_{kN-|\alpha|} 
\otimes v_{kN-|\alpha|+1}\dots v_{kN}\beta) \\ 
&= 
1\otimes \alpha v_1\dots v_{kN-|\alpha|} 
\otimes v_{kN-|\alpha|+1}\dots v_{kN}\beta. 
\end{align*} 
Hence 
$d_{2k+1}r_{2k}+r_{2k-1}d_{2k}(\alpha\otimes v_1\dots v_{kN}\otimes\beta) 
=\alpha\otimes v_1\dots v_{kN}\otimes \beta. 
$ 
 
\vspace{2mm} 
 
Similarly, for $n=2k+1$ we have 
$d_{2k+2}r_{2k+1}(\alpha\otimes v_1\dots v_{kN+1}\otimes\beta)=0$ if $|\alpha|<N-1$ and 
\begin{align*} 
d_{2k+2}&r_{2k+1}(\alpha\otimes v_1\dots v_{kN+1}\otimes\beta)\\ 
&= 
d_{2k+2}(1\otimes \alpha v_1\dots v_{kN+1}\otimes\beta)\\ 
&= 
\sum_{j=0}^{N-1} 
a_1\dots a_{j}\otimes 
\underbrace{a_{j+1}\dots a_{N-1}v_1\dots v_{t}}_{kN+1} 
\otimes v_{t+1}\dots v_{kN+1}\beta, 
\end{align*} 
if $|\alpha|=N-1$. 
On the other hand 
\begin{multline*} 
r_{2k}d_{2k+1}(\alpha\otimes v_1\dots v_{kN+1}\otimes\beta)\\ 
= 
r_{2k}(\alpha v_1\otimes v_2\dots v_{kN+1}\otimes\beta 
-\alpha\otimes v_1\dots v_{kN}\otimes v_{kN+1}\beta). 
\end{multline*} 
If $|\alpha|<N-1$ then 
$r_{2k}(\alpha v_1\otimes v_2\dots v_{kN+1}\otimes\beta 
-\alpha\otimes v_1\dots v_{kN}\otimes v_{kN+1}\beta)$ 
is a telescopic sum that adds up to $\alpha\otimes v_1\dots v_{kN+1}\otimes \beta$. 
 
If $|\alpha|=N-1$ then $\alpha v_1\otimes v_2\dots v_{kN+1}\otimes \beta=0$ and 
\begin{align*} 
r_{2k}d_{2k+1}(\alpha\otimes v_1&\dots v_{kN+1}\otimes \beta)\\ 
&=- 
r_{2k}(\alpha\otimes v_1\dots v_{kN}\otimes v_{kN+1}\beta) \\ 
&=- 
\sum_{j=1}^{N-1} 
a_1\dots a_{j-1}\otimes 
\underbrace{a_{j}\dots a_{N-1}v_1\dots v_{t}}_{kN+1} 
\otimes v_{t+1}\dots v_{kN+1}\beta. 
\end{align*} 
Hence 
$d_{2k+2}r_{2k+1}+r_{2k}d_{2k+1}(\alpha\otimes v_1\dots v_{kN+1}\otimes\beta) 
=\alpha\otimes v_1\dots v_{kN+1}\otimes \beta. 
$ 
\end{proof} 
 
\ 
 
\section{The comparison morphisms} 
\label{sec:comparison-morphisms} 
 
A comparison morphism between two projective resolutions of an algebra $A$ 
is a morphism of chain complexes that lifts the identity map on $A$. 
Such a morphism induces a quasi-isomorphism between the derived complexes 
$\Hom_{A^e}(\cdot,A)$. 
 
In this section we define maps 
\[ 
 {\bf F}:{\bf P}\to{\bf Q}\quad \text{ and }\quad{\bf G}:{\bf Q}\to{\bf P} 
\] 
between these $A$-bimodule resolutions of $A$ and we state in 
Theorem \ref{MainComparison} that they are in fact comparison morphisms. 
This is one of the main results of the paper. 
The proofs, which we find non trivial and subtle, are in 
Sections \ref{Comp G} and \ref{Comp F}. 
The reader interested only on the main results in the paper may safely skip 
these two sections. 
 
We define ${\bf F}$ and ${\bf G}$ as the $A$-bimodule extensions of 
maps defined on elements of $\mathcal B'_{{\bf P}_{n}}$ and 
$\mathcal B'_{{\bf Q}_{n}}$ respectively. 
As in the case of the differentials $b$ and $d$ one should check that these 
$A$-bimodule extensions are well defined, but this is straightforward since 
the tensor products in ${\bf Q}$ and ${\bf P}$ are both over 
${\mathrm{k}}\Delta_0$.

\subsection{The comparison morphism ${\bf F}:{\bf P}\to{\bf Q}$}\label{subsec:F} 
Let ${\bf F}_{0}=id$ and, for $n\ge1$, let ${\bf F}_n:{\bf P}\to{\bf Q}$ be 
the $A$-bimodule extension of the following map defined on elements of 
$\mathcal{B}'_{{\bf P}_{n}}$. 
If $n=2k$ and $p=1\otimes v_1\dots v_{kN}\otimes1\in\mathcal B'_{{\bf P}_{2k}}$, 
let 
\begin{multline*} 
{\bf F}_{2k}(p)= \\ 
\sum 
1 
[ 
\underbrace{v_1\dots v_{x_1}}_{x_1} 
| 
\underbrace{v_{1+x_1}}_{1} 
| 
\underbrace{\dots v_{1+x_1+x_2}}_{x_2} 
| 
\underbrace{v_{2+x_1+x_2}}_{1} 
|\dots\dots| 
\underbrace{v_{k+\sum x_j}}_{1} 
] 
\underbrace{\dots\dots v_{kN}}_{kN-k-\sum x_j}, 
\end{multline*} 
where the sum is taken over all $k$-tuples $(x_1,\dots,x_k)\in\Z^k$ 
such that $1\le x_i<N$. 
If $n=2k+1$ and $p=1\otimes v_1\dots v_{kN+1}\otimes1\in\mathcal B'_{{\bf P}_{2k+1}}$, 
let 
\begin{multline*} 
{\bf F}_{2k+1}(p)= \\ 
\sum 
1 
[\! 
\underbrace{\!\!v_1\!\!}_{1}\! 
| 
\underbrace{v_2\dots v_{1+x_1}}_{x_1} 
| 
\underbrace{v_{2+x_1}}_{1} 
| 
\underbrace{\dots v_{2+x_1+x_2}}_{x_2} 
|\dots\dots| 
\underbrace{v_{k\!+\!1+\sum x_j}}_{1} 
] 
\underbrace{\dots\dots v_{kN+1}}_{kN-k-1-\sum x_j}, 
\end{multline*} 
where the sum is taken over all $k$-tuples $(x_1,\dots,x_k)\in\Z^k$ 
such that $1\le x_i<N$. 
 
In order to prove that ${\bf F}$ is a comparison morphism, one should first 
check that the diagram 
\[ 
\xymatrix{ 
A\otimes_{\Delta_{0}}A_+\otimes_{\Delta_{0}}A\ar[r]^{\qquad b_{1}}& 
A\otimes A\ar[d]^{\text{id}}\ar[r]^{\epsilon}& 
A 
                                        \\ 
A\otimes_{\Delta_{0}}{\mathrm{k}}\Delta_{1}\otimes_{\Delta_{0}}A\ar[u]^{{\bf F}_{1}}\ar[r]_{\qquad d_1}& 
A\otimes A\ar[r]^{\epsilon}&A. 
} 
\] 
is commutative. This is immediate since ${\bf F}_{1}(1\otimes v\otimes1)=1[v]1$.

\subsection{The comparison morphism ${\bf G}:{\bf Q}\to{\bf P}$} 
\label{subsec:G} 
 
Let ${\bf G}_{0}=id$ and, for $n\ge1$, let ${\bf G}_n:{\bf Q}\to{\bf P}$ be 
the $A$-bimodule extension of the following map defined on elements of 
$\mathcal{B}'_{{\bf Q}_{n}}$. 
For 
$q=1[\alpha_1|\dots|\alpha_n]1\in\mathcal B'_{{\bf Q}_n}$, 
let $v(q)$ be the arrow decomposition of the path $\alpha_1\alpha_2\dots \alpha_n$ in 
${\mathrm{k}}\Delta$ (not in $A$). 
Thus, if $\alpha_i=a^i_1\dots a_{|\alpha_i|}^i$, then 
\begin{equation*} 
v(q)=v_1\dots v_{|q|}=a_1^1\dots a_{|\alpha_1|}^1a_1^2\dots a_{|\alpha_2|}^2\dots\dots a_1^{n} 
\dots a_{|\alpha_{n}|}^{n}\in\Delta_{|q|}, 
\end{equation*} 
where $|q|=\sum_{i=1}^{n}|\alpha_i|$. 
Note that $v(q)\ne0$ for all $q\in\mathcal B'_{{\bf Q}_n}$. 
 
If $n=2k$ and $v(q)=v_1\dots v_{|q|}$, 
let 
\begin{align*} 
{\bf G}_{2k}(q)= 
\begin{cases} 
1\otimes v_1\dots v_{kN}\otimes v_{kN+1}\dots v_{|q|},& 
\text{if }\alpha_{2i-1}\alpha_{2i}=0\text{ for }i=1\dots k; \\ 
0, &\text{otherwise}. 
\end{cases} 
\end{align*} 
Note that the condition $\alpha_{2i-1}\alpha_{2i}=0$ for all $i=1\dots k$ implies 
that $|q|\ge kN$. 
 
Similarly, if $n=2k+1$ and 
$v(q)=v_1\dots v_{|q|}$, 
let 
\begin{align*} 
{\bf G}_{2k+1}(q)= 
\begin{cases} 
\displaystyle\sum_{j=1}^{|\alpha_1|} 
v_1\dots v_{j-1}\otimes v_{j}\dots v_{kN+j}\otimes v_{kN+j+1}\dots v_{|q|},& 
\text{if }\alpha_{2i}\alpha_{2i+1}=0\\[-5mm] 
& \text{for }i=1\dots k; \\[2mm] 
0, &\text{otherwise}. 
\end{cases} 
\end{align*} 
Since $|\alpha_1|\ge1$, then $|q|\ge kN+1$ provided that $\alpha_{2i}\alpha_{2i+1}=0$ 
for all $i=1\dots k$. 
Note also that only $\alpha_1$ is involved in the sum. 
 
In this case, the commutativity of the first diagram 
\[ 
\xymatrix{ 
A\otimes_{\Delta_{0}}A_+\otimes_{\Delta_{0}}A 
\ar[d]_{{\bf G}_{1}}\ar[r]^{\qquad b_{1}}& 
A\otimes A\ar[d]^{\text{id}}\ar[r]^{\epsilon}& 
A 
                                        \\ 
A\otimes_{\Delta_{0}}{\mathrm{k}}\Delta_{1}\otimes_{\Delta_{0}}A\ar[r]_{\qquad d_1}& 
A\otimes A\ar[r]^{\epsilon}&A 
} 
\] 
follows immediately by evaluating the telescopic sum $d_1\circ{\bf G}_1$.

\begin{theorem}\label{MainComparison} 
The following diagram is commutative for all $k\ge1$ and therefore 
${\bf F}$ and ${\bf G}$ are comparison morphisms between the $A^e$-projective resolutions 
${\bf P}$ and ${\bf Q}$. 
 
\[ 
\xymatrix 
{ 
\scriptstyle 
A\otimes_{{}_{\!\Delta_{0}}}\!\! 
A_+^{\otimes^{2k+1}_{\Delta_0}} 
\otimes_{{}_{\!\Delta_{0}}}\!\!\!\!A 
& 
\scriptstyle 
A\otimes_{{}_{\!\Delta_{0}}}\!\! 
A_+^{\otimes^{2k}_{\Delta_0}} 
\otimes_{{}_{\!\Delta_{0}}}\!\!\!\!A 
& 
\scriptstyle 
A\otimes_{{}_{\!\Delta_{0}}}\!\! 
A_+^{\otimes^{2k-1}_{\Delta_0}} 
\otimes_{{}_{\!\Delta_{0}}}\!\!\!\!A 
                                       \\ 
{\bf Q}_{2k+1}\ar@{}|\|[u]\ar[r]^{b_{2k+1}}\ar[d]<2pt>^{{\bf G}_{2k+1}}& 
{\bf Q}_{2k}\ar@{}|\|[u]\ar[r]^{b_{2k}}\ar[d]<2pt>^{{\bf G}_{2k}}& 
{\bf Q}_{2k-1}\ar@{}|\|[u]\ar[d]<2pt>^{{\bf G}_{2k-1}} 
                                       \\ 
{\bf P}_{2k+1}\ar[r]_{d_{2k+1}}\ar[u]<2pt>^{{\bf F}_{2k+1}}& 
{\bf P}_{2k}\ar[r]_{d_{2k}}\ar[u]<2pt>^{{\bf F}_{2k}}& 
{\bf P}_{2k-1}\ar[u]<2pt>^{{\bf F}_{2k-1}} 
                                        \\ 
\scriptstyle 
A\otimes_{{}_{\!\Delta_{0}}}\!\!\!\!{\mathrm{k}}\Delta_{k\!N\!+\!1}\otimes_{{}_{\!\Delta_{0}}}\!\!\!\!A 
\ar@{}|\|[u] 
& 
\scriptstyle 
A\otimes_{{}_{\!\Delta_{0}}}\!\!\!\!{\mathrm{k}}\Delta_{k\!N}\otimes_{{}_{\!\Delta_{0}}}\!\!\!\!A 
\ar@{}|\|[u] 
& 
\scriptstyle 
A\otimes_{{}_{\!\Delta_{0}}}\!\!\!\!{\mathrm{k}}\Delta_{(k\!-\!1)\!N\!+\!1}\otimes_{{}_{\!\Delta_{0}}}\!\!\!\!A 
\ar@{}|\|[u] 
} 
\] 
\end{theorem}

\section{${\bf G}:{\bf Q}\to{\bf P}$ is a comparison morphism} 
\label{sec:proof-G} 
 
The proof is divided into two parts, (A) and (B), 
corresponding to the cases $n$ even and $n$ odd respectively. 
 
(A) Assume $n=2k$. Let $q=1[\alpha_1|\dots|\alpha_{2k}]1\in\mathcal B'_{{\bf Q}_{2k}}$ and let 
$v(q)=v_1\dots v_{|q|}$. 
Let 
\[ 
M=\{i\in\{1,\dots,k\}:\alpha_{2i-1}\alpha_{2i}\ne0\}. 
\] 
 
\vskip.4cm 
 
\noindent 
(A1) Case $M=\emptyset$. 
We have 
\[ 
{\bf G}_{2k}(q)= 
1\otimes v_1\dots v_{kN}\otimes v_{kN+1}\dots v_{|q|} 
\] 
and 
\[ 
d_{2k}({\bf G}_{2k}(q))= 
\sum_{j=0}^{N-1} 
v_1\dots v_j\otimes 
\underbrace{v_{j+1}\dots v_{t}}_{(k-1)N+1}\otimes v_{t+1}\dots v_{|q|}. 
\] 
On the other hand, $M=\emptyset$ implies 
\[ 
b_{2k}(q)= 
\alpha_1[\alpha_2|\dots|\alpha_{2k}]1+ 
\sum_{i=1}^{k-1}1[\alpha_1|\dots|\alpha_{2i}\alpha_{2i+1}|\dots| \alpha_{2k}]1+ 
1[\alpha_1|\dots|\alpha_{2k-1}]\alpha_{2k} 
\] 
and, since $N\le |\alpha_1|+|\alpha_2|$, 
\begin{multline*} 
{\bf G}_{2k-1}(b_{2k}(q)) = 
\sum_{j=|\alpha_1|+1}^{N} 
v_1\dots v_{j-1}\otimes v_{j}\dots v_{(k-1)N+j}\otimes v_{(k-1)N+1+j}\dots v_{|q|}\\ 
+ 
\sum_{i=1}^{k-1}{\bf G}_{2k-1}(1[\alpha_1|\dots|\alpha_{2i}\alpha_{2i+1}|\dots| \alpha_{2k}]1)+ 
{\bf G}_{2k-1}(1[\alpha_1|\dots|\alpha_{2k-1}]\alpha_{2k}). 
\end{multline*} 
In the second line, all  terms but one are zero, depending on which is the first 
$j$, $1\le j\le k-1$, for which $\alpha_{2j}\alpha_{2j+1}=0$. 
The nonzero term is 
\[ 
\sum_{j=1}^{|\alpha_1|} 
v_1\dots v_{j-1}\otimes v_{j}\dots v_{(k-1)N+j}\otimes v_{(k-1)N+1+j}\dots v_{|q|} 
\] 
and therefore 
${\bf G}_{2k-1}(b_{2k}(q))=d_{2k}({\bf G}_{2k}(q))$. 
 
\vskip.4cm 
 
\noindent 
(A2) Case  $M\ne\emptyset$. We have now $d_{2k}({\bf G}_{2k}(q))=0$. 
Let 
\[ 
M'=\{i\in\{1,\dots,k\}:{\bf G}_{2k-1}(1[\alpha_1|\dots|\alpha_{2i-1}\alpha_{2i}|\dots |\alpha_{2k}]1)\ne0\} 
\subset M. 
\] 
 
\vskip.2cm 
 
\noindent 
(A2a) Assume $M'\ne\emptyset$ and let $i_0$ be the smallest element in $M'$. 
We assume that $i_0>1$ since the case $i_0=1$ is easier. 
By the definition of ${\bf G}_{2k-1}$, it follows that 
\begin{align*} 
\alpha_{2i}\alpha_{2i+1}&=0,\text{ for }i=1\dots i_0-2, \text{ and }\\ 
\alpha_{2i-1}\alpha_{2i}&=0,\text{ for }i=i_0+1\dots k. 
\end{align*} 
In particular $i_0$ is the largest element of $M$ and 
\begin{multline*} 
{\bf G}_{2k-1}(b_{2k}(q)) = \\ 
{\bf G}_{2k-1}(\alpha_1[\alpha_2|\dots|\alpha_{2k}]1)-{\bf G}_{2k-1}(1[\alpha_1|\alpha_2|\dots|\alpha_{2i_0-1}\alpha_{2i_0}|\dots |\alpha_{2k}]1)\\ 
+ 
\sum_{i=i_0-1}^{k-1}{\bf G}_{2k-1}(1[\alpha_1|\dots|\alpha_{2i}\alpha_{2i+1}|\dots| \alpha_{2k}]1)+ 
{\bf G}_{2k-1}(1[\alpha_1|\dots|\alpha_{2k-1}]\alpha_{2k}); 
\end{multline*} 
As in the case (A1), all terms but one are zero in the second line, 
and this term cancels out with the first line.

\vskip.2cm 
 
\noindent 
(A2b) Assume $M'=\emptyset$. Thus ${\bf G}_{2k-1}(b_{2k}(q))$ contains only positive terms and 
we must prove that all of them are zero. 
Let $i_0$ be any element of $M$. Since $\alpha_{2i_0-1}\alpha_{2i_0}\ne0$ the definition of 
${\bf G}_{2k-1}$ implies that 
\[ 
{\bf G}_{2k-1}(1[\alpha_1|\dots|\alpha_{2i}\alpha_{2i+1}|\dots| \alpha_{2k}]1)=0 
\text{ for all }i=1\dots i_0-1. 
\] 
We now take $i_0=\max(M)$. 
Since $i_0\not\in M'$, the maximality of $i_0$ implies that either 
$\alpha_{2(i_0-1)}\alpha_{2i_0-1}\alpha_{2i_0}\ne0$ or there exists 
$j_0<i_0-1$ such that $\alpha_{2j_0}\alpha_{2j_0+1}\ne0$. 
In any case, there exists $j_0<i_0$ such that $\alpha_{2j_0}\alpha_{2j_0+1}\ne0$ (in particular $i_0>1$). 
Therefore 
${\bf G}_{2k-1}(1[\alpha_1|\dots|\alpha_{2i}\alpha_{2i+1}|\dots| \alpha_{2k}]1)=0$ for all $i=j_0+1\dots k-1$. 
Since the extreme cases 
${\bf G}_{2k-1}(\alpha_1[\alpha_2|\dots| \alpha_{2k}]1)$ and ${\bf G}_{2k-1}(1[\alpha_1|\dots| \alpha_{2k-1}]\alpha_{2k})$ 
are clearly zero too, this completes the proof in case (A).

\vskip.4cm 
 
\noindent 
(B) Assume $n=2k+1$. Let $q=1[\alpha_1|\dots|\alpha_{2k+1}]1\in\mathcal B'_{{\bf Q}_{2k}}$ and let 
$v(q)=v_1\dots v_{|q|}$. 
Let 
\[ 
M=\{i\in\{1,\dots,k\}:\alpha_{2i}\alpha_{2i+1}\ne0\}. 
\] 
 
\vskip.4cm 
 
\noindent 
(B1) Case $M=\emptyset$. 
Then 
\[ 
{\bf G}_{2k+1}(q)= 
\sum_{j=1}^{|\alpha_1|} 
v_1\dots v_{j-1}\otimes v_{j}\dots v_{kN+j}\otimes v_{kN+1+j}\dots v_{|q|} 
\] 
and 
\begin{multline*} 
d_{2k+1}({\bf G}_{2k+1}(q))= 
v_1\dots v_{|\alpha_1|}\otimes 
v_{|\alpha_1|+1}\dots v_{|\alpha_1|+kN}\otimes v_{|\alpha_1|+kN+1}\dots v_{|q|} \\ 
-1\otimes v_{1}\dots v_{kN}\otimes v_{kN+1}\dots v_{|q|}. 
\end{multline*} 
On the other hand, $M=\emptyset$ implies 
\begin{multline*} 
b_{2k+1}(q)= 
\alpha_1[\alpha_2|\dots|\alpha_{2k+1}]1 \\ 
- \sum_{i=1}^{k}1[\alpha_1|\dots|\alpha_{2i-1}\alpha_{2i}|\dots| \alpha_{2k+1}]1 
-1[\alpha_1|\dots|\alpha_{2k}]\alpha_{2k+1} 
\end{multline*} 
and, 
\begin{multline*} 
{\bf G}_{2k}(b_{2k+1}(q)) = 
v_1\dots v_{|\alpha_1|}\otimes 
v_{|\alpha_1|+1}\dots v_{|\alpha_1|+kN}\otimes v_{|\alpha_1|+kN+1}\dots v_{|q|} \\ 
- \sum_{i=1}^{k}{\bf G}_{2k}(1[\alpha_1|\dots|\alpha_{2i-1}\alpha_{2i}|\dots| \alpha_{2k+1}]1) 
-{\bf G}_{2k}(1[\alpha_1|\dots|\alpha_{2k}]\alpha_{2k+1}). 
\end{multline*} 
As in case (A),  all  terms but one are zero in the second line, 
and the nonzero term is 
$1\otimes v_{1}\dots v_{kN}\otimes v_{kN+1}\dots v_{|q|}$. 
This yields 
${\bf G}_{2k}(b_{2k+1}(q))=d_{2k+1}({\bf G}_{2k+1}(q))$. 
 
\vskip.4cm 
 
\noindent 
(B2) Case $M\ne\emptyset$. We have $d_{2k}({\bf G}_{2k}(q))=0$. 
Let 
\[ 
M'=\{i\in\{1,\dots,k\}:{\bf G}_{2k}(1[\alpha_1|\dots|\alpha_{2i}\alpha_{2i+1}|\dots |\alpha_{2k+1}]1)\ne0\} 
\subset M. 
\] 
 
\vskip.2cm 
 
\noindent 
(B2a) Assume $M'\ne\emptyset$ and let $i_0$ be the largest element in $M'$. 
By the definition of ${\bf G}_{2k}$, it follows that 
\begin{align*} 
\alpha_{2i-1}\alpha_{2i}&=0,\text{ for }i=1\dots i_0-1, \text{ and }\\ 
\alpha_{2i}\alpha_{2i+1}&=0,\text{ for }i=i_0+1\dots k. 
\end{align*} 
In particular $i_0$ is the smallest element of $M$ and 
\begin{multline*} 
{\bf G}_{2k+1}(b_{2k}(q)) = 
{\bf G}_{2k}(1[\alpha_1|\alpha_2|\dots|\alpha_{2i_0}\alpha_{2i_0+1}|\dots |\alpha_{2k+1}]1)\\ 
- 
\sum_{i=i_1}^{k-1}{\bf G}_{2k}(1[\alpha_1|\dots|\alpha_{2i-1}\alpha_{2i}|\dots| \alpha_{2k+1}]1)- 
{\bf G}_{2k}(1[\alpha_1|\dots|\alpha_{2k}]\alpha_{2k+1}); 
\end{multline*} 
where $i_1=i_0-1$, if $i_0>1$; and $i_1=1$, if $i_0=1$. 
All terms but one are zero in the second line, 
and this term cancels out with the first line.

\vskip.2cm 
 
\noindent 
(B2b) Assume $M'=\emptyset$. Thus ${\bf G}_{2k}(b_{2k+1}(q))$ contains only negative terms and 
we must prove that all of them are zero. 
Let $i_0$ be any element of $M$. Then 
\[ 
{\bf G}_{2k}(1[\alpha_1|\dots|\alpha_{2i-1}\alpha_{2i}|\dots| \alpha_{2k+1}]1)=0 
\text{ for all }i=1\dots i_0-1. 
\] 
We now take $i_0=\max(M)$. 
Since $i_0\not\in M'$, the maximality of $i_0$ implies that either 
$\alpha_{2i_0-1}\alpha_{2i_0}\alpha_{2i_0+1}\ne0$ or there exist $j_0<i_0$ such that 
$\alpha_{2j_0-1}\alpha_{2j_0}\ne0$. 
In any case,  there exist $j_0\le i_0$ such that $\alpha_{2j_0-1}\alpha_{2j_0}\ne0$. 
Therefore 
${\bf G}_{2k}(1[\alpha_1|\dots|\alpha_{2i-1}\alpha_{2i}|\dots| \alpha_{2k+1}]1)=0$ for all $i=j_0\dots k$. 
Since the extreme case 
${\bf G}_{2k}(1[\alpha_1|\dots| \alpha_{2k}]\alpha_{2k})$ 
is clearly zero too, this completes the case (B) and the proof.

\section{${\bf F}:{\bf P}\to{\bf Q}$ is a comparison morphism} 
\label{sec:proof-F} 
 
We need some preliminary results. 
 
\subsection{A complex of compositions} 
 
Let $C_n(m)$ be the set of all the compositions (ordered partitions) of $m$ in $n$ parts in 
which only the first and last parts are allowed to be zero. 
That is, 
\[ 
C_n(m)=\{[c_1,\dots,c_n]:c_j\in\N_0, 
\,c_j>0\text{ for $j=2\dots n-1$ and }\Sigma c_i=m\}. 
\] 
Let ${\bf C}_n(m)={\mathrm{k}} C_n(m)$ be the  vector space with basis $C_n(m)$. 
 
If ${\alpha}=[x_1,\dots,x_n]\in C_n(m)$ and 
${\beta}=[y_1,\dots,y_{n'}]\in C_{n'}(m')$ 
then we shall denote by 
$[\alpha,\beta]$ the juxtaposition of $\alpha$ and $\beta$, 
that is 
\[ 
[\alpha,\beta]= 
[x_1,\dots,x_n,y_1,\dots,y_{n'}] 
\in C_{n+n'}(m+m'). 
\] 
Analogously, if ${\alpha}\in{C}_{n}(m)$ and 
${\beta}\in{\bf C}_{n'}(m')$, with $\beta=\sum \beta_j$, 
${\beta_i}\in C_{n'}(m')$, then 
$[\alpha,\beta]=\sum[\alpha,\beta_j]$. 
 
Let  ${D}:{\bf C}_n(m)\to{\bf C}_{n-1}(m)$ be the usual differential of 
compositions, 
\[ 
{D}([x_1,\dots,x_n])=\sum_{j=1}^{n-1}(-1)^{j+1} 
[x_1,\dots,x_j+x_{j+1},\dots,x_n]. 
\] 
It is straightforward to see that ${D}^2=0$. 
Moreover, if ${\bf W}_n(m,N)\subset{\bf C}_n(m)$ is the subspace spanned by 
the compositions containing some part larger than or equal to $N$, 
then ${D}({\bf W}_n(m,N))\subset {\bf W}_n(m,N)$ 
and thus ${D}$ factors through  the quotient 
\[ 
{\bf C}_n(m,N)={\bf C}_n(m)/{\bf W}_n(m,N). 
\] 
For ${\alpha}\in{\bf C}_n(m)$, let 
${c_N}\in{\bf C}_n(m,N)$ be its projection and let, by definition, 
${D}_{N}({c_N})={D}({c})_{N}$. 
Thus $({\bf C}_n(m,N),{D}_N)$ is again a complex. 
 
\ 
 
Let $I_N=\{1,2,\dots,N-1\}$ and 
for  $k\ge1$ and $M\ge k(N-1)$ 
let 
\begin{align*} 
\alpha_{M}^{k}&:I_N^k\to {\bf C}_{2k+1}({M\!+\!k},N),\; 
&&\alpha_{M}^{k}(x)=[x_1,1,x_2,1,\dots,x_k,1,M-\Sigma x_i]_{N} \\ 
\beta_{M}^{k}&:I_N^k\to {\bf C}_{2k+2}({M\!+\!k+\!1},N),\; 
&&\beta_{M}^{k}({x})=[1,x_1,1,x_2,1,\dots,x_k,1 ,M-\Sigma x_i]_{N}. 
\end{align*} 
The assumption $M\ge k(N-1)$ is necessary in order to ensure a non negative 
last part for any ${x}\in I_N^k$. 
 
For $k\ge1$ let 
\begin{align*} 
{A_{M}^{k}}& 
=\sum_{{ x}\in I_N^k}\alpha_M^{k}({ x})\in{\bf C}_{2k+1}(M\!+\!k,N), \\ 
{B_{M}^{k}}& 
=\sum_{{ x}\in I_N^k}\beta_M^{k}({ x})\in{\bf C}_{2k+2}(M\!+\!k\!+1,N). 
\end{align*} 
We define $B_{M}^0=[1,M]_{N}$ and we need not define $A_{M}^k$ for $k=0$. 
 
\begin{lemma}\label{lemma:comp1} For all $k\ge1$ and $M\ge k(N-1)$ we have 
\[ 
{B_{M}^{k}}=[1,A^{k}_{M}]_{N}\; 
\text{ and }\; 
{A_{M}^{k}}=\sum_{j=1}^{N-1}[j,B^{k-1}_{M-j}]_{N}. 
\] 
Moreover, ${A_{M}^{k}}\ne0$ (resp. ${B_{M}^{k}}\ne0$) 
if and only if $M\le (k+1)(N-1)$. 
\end{lemma} 
\begin{proof} 
The first part of the lemma is straightforward from the definition of 
$\alpha_{M}^{k}({ x})$ and $\beta_{M}^{k}({ x})$. 
The second part follows from the fact that 
the last part of the composition $\alpha_{M}^{k}({x})$ 
(resp. $\beta_{M}^{k}({ x})$) is greater than or equal to $N$ for all 
${ x}\in I_N^k$, if and only if $M>(k+1)(N-1)$. 
\end{proof} 
 
\begin{lemma}\label{lemma:comp2} 
For all $k\ge1$ and $M\ge k(N-1)$ we have 
\[ 
{D_N}(A_M^k)=-B_{M}^{k-1}\quad\text{ and }\quad{D_N}(B_M^k)=A_{M+1}^k. 
\] 
\end{lemma} 
 
\noindent 
{\bf Remark.} 
The following picture shows the values of $M$ for which $A$ and $B$ are 
different from zero. In particular it shows that the lemma is 
consistent with the fact that ${D}^2_{N}=0$.

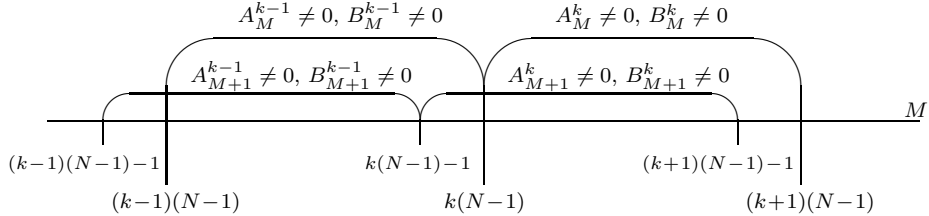
\begin{center} 
\setlength{\unitlength}{30pt} 
\begin{picture}(4,4)(-2,-2) 
  \put(-5.5,0){\line(1,0){11}} 
  \put(0,-.8){\line(0,1){1.2}} 
  \put(-0.5,-1.1){\SMALL $k(N\!-\!1)$} 
  \put(-4,-.8){\line(0,1){1.2}} 
  \put(-4.7,-1.1){\SMALL $(k\!-\!1)(N\!-\!1)$} 
  \put(4,-.8){\line(0,1){1.2}} 
  \put(3.3,-1.1){\SMALL $(k\!+\!1)(N\!-\!1)$} 
  \put(-4.8,-.3){\line(0,1){.3}} 
  \put(-6,-.6){\tiny $(k\!-\!1)(N\!-\!1)\!-\!1$} 
  \put(-.8,-.3){\line(0,1){.3}} 
  \put(-1.5,-.6){\tiny $k(N\!-\!1)\!-\!1$} 
  \put(3.2,-.3){\line(0,1){.3}} 
  \put(2,-.6){\tiny $(k\!+\!1)(N\!-\!1)\!-\!1$} 
  \put(2,.4){\oval(4,1.3)[t]} 
  \put(-2,.4){\oval(4,1.3)[t]} 
  \put(-2.8,0){\oval(4,.7)[t]} 
  \put(1.2,0){\oval(4,.7)[t]} 
  \put(5.3,0.05){\SMALL $M$} 
  \put(.9,1.25){\SMALL $A_M^k\ne0$, $B_M^k\ne0$} 
  \put(-3.1,1.25){\SMALL $A_M^{k-1}\ne0$, $B_M^{k-1}\ne0$} 
  \put(.3,.5){\SMALL $A_{M+1}^k\ne0$, $B_{M+1}^k\ne0$} 
  \put(-3.7,.5){\SMALL $A_{M+1}^{k-1}\ne0$, $B_{M+1}^{k-1}\ne0$} 
\end{picture} 
\end{center}

\begin{proof} 
We assume that $M$ is fixed. 
We shall now prove simultaneously both equalities by induction on $k$ for 
$1\le k\le\frac{M}{N-1}$. 
 
If $k=1$ then 
\begin{align*} 
{D_N}(A_{M}^1)&= 
\sum_{x_1=1}^{N-1} [x_1+1,M-x_1]_{N}- 
\sum_{x_1=1}^{N-1} [x_1,M-x_1+1]_{N}\\ 
&=- [1,M]_{N} \\ 
&=-B_{M}^0. 
\end{align*} 
This completes the case $k=1$ for the first equation. 
The case $k=1$ for the second equation is analogous: 
\[ 
B_M^1= 
\sum_{x_1=1}^{N-1} [1,x_1,1,M-x_1]_{N} 
\] 
and 
\begin{align*} 
{D_N}(B_M^1)= 
\sum_{x_1=1}^{N-1} [x_1+1,1,M-x_1]_{N}- 
& 
\sum_{x_1=1}^{N-1} [1,x_1+1,M-x_1]_{N} \\ 
& 
+\sum_{x_1=1}^{N-1} [1,x_1,M-x_1+1]_{N}. 
\end{align*} 
The last two terms add up as a telescopic sum and the result is 
\[ 
[1,1,M]_{N}-[1,N,M-N+1]_{N}=[1,1,M]_{N}. 
\] 
Therefore 
\begin{align*} 
{D_N}(B_M^1) 
&= 
\sum_{x_1=1}^{N-1} [x_1+1,1,M-x_1]_{N}+[1,1,M]_{N} \\ 
&= 
\sum_{{ x}\in I_N^1}\alpha_{M+1}^1({ x}) \\ 
&=A_{M+1}^1. 
\end{align*} 
 
Now we assume that the lemma is true for $k-1$. 
From Lemma \ref{lemma:comp1} we have 
\[ 
{A_{M}^{k}} 
=\sum_{j=1}^{N-1} [j,B_{M-j}^{k-1}]_{N} 
=\sum_{j=1}^{N-1} [j,1,A_{M-j}^{k-1}]_{N}. 
\] 
Hence 
\begin{align*} 
{D_N}(A_M^k) 
&= 
\sum_{j=1}^{N-1} 
[1+j,A_{M-j}^{k-1}]_{N}- 
\sum_{j=1}^{N-1} 
[j,{D_N}(B^{k-1}_{M-j})]_{N} 
                                     \\ 
&= 
\sum_{j=2}^{N-1} 
[j,A_{M+1-j}^{k-1}]_{N}- 
\sum_{j=1}^{N-1}[j,A^{k-1}_{M+1-j}]_{N} 
                                     \\ 
&= -[1,A_{M}^{k-1}]_{N} 
                                       \\ 
&= 
-B_{M}^{k-1}. 
\end{align*} 
 
Similarly, 
\[ 
{B_{M}^{k}} 
=[1,A^{k}_{M}]_{N} 
=\sum_{j=1}^{N-1} [1,j,1,A_{M-j}^{k-1}]_{N}. 
\] 
 
Hence 
\begin{align*} 
{D_N}(B_M^k) 
&= 
\sum_{j=1}^{N-1} 
[1+j,1,A_{M-j}^{k-1}]_{N}- 
[1,{D_N}(A^{k}_{M})]_{N} 
                                     \\ 
&= 
\sum_{j=1}^{N-2} 
[1+j,1,A_{M-j}^{k-1}]_{N}- 
[1,-B^{k-1}_{M}]_{N} 
                                     \\ 
&= 
\sum_{j=1}^{N-1} 
[j,B_{M+1-j}^{k-1}]_{N} 
                                       \\ 
&= 
A_{M+1}^k. 
\end{align*} 
This completes the inductive argument. 
\end{proof}

Now we define 
\begin{align*} 
\tilde A_{M}^k&=[0,A_M^k]_{N}\in{\bf C}_{2k+2}(M\!+\!k,N), 
&&\text{for }k\ge1\text{ and }M\ge k(N-1), \\ 
\tilde B_{M}^k&=[0,B_M^k]_{N}\in{\bf C}_{2k+3}(M\!+\!k+\!1,N), 
&&\text{for }k\ge0\text{ and }M\ge k(N-1). \\ 
\end{align*} 
 
\begin{proposition}\label{Prop D} 
Let $k\ge1$ and $M\ge k(N-1)$ we have 
\begin{align*} 
{D_N}(\tilde B_M^k)&=[1,A_M^k]_{N}-[0,A_M^k]_{N}\text{ and }\\ 
{D_N}(\tilde A_M^k)&=\sum_{j=0}^{N-1}[j,B_{M-j}^{k-1}]_{N}. 
\end{align*} 
\end{proposition} 
 
\begin{proof} 
Since $\tilde B_M^k=[0,B^k_M]_{N}=[0,1,A^k_M]_{N}$, 
then Lemma \ref{lemma:comp2} implies 
\begin{align*} 
{D_N}(\tilde B_M^k) 
&=[1,A^k_M]_{N}-[0,{D_N}(B^k_M)]_{N} \\ 
&=[1,A^k_M]_{N}-[0,A^k_M]_{N}. 
\end{align*} 
Similarly, since 
$\tilde A_M^k=[0,A^k_M]_{N}=\sum_{j=1}^{N-1}[0,j,B^{k-1}_{M-j}]_{N}$, 
then 
\begin{align*} 
{D_N}(\tilde A_M^k) 
&=\sum_{j=1}^{N-1}[j,B^{k-1}_{M-j}]_{N}-[0,{D_N}(A^k_M)]_{N} \\ 
&=\sum_{j=0}^{N-1}[j,B^{k-1}_{M-j}]_{N}. 
\end{align*} 
\end{proof} 
 
\subsection{The final step} 
 
Each composition 
${\alpha}=[x_1,\dots,x_n]\in C_n(m)$ defines an $A$-bimodule morphism 
\[ 
\phi_{\alpha}: 
A\otimes_{{\mathrm{k}}\Delta_{0}}{\mathrm{k}}\Delta_{m}\otimes_{{\mathrm{k}}\Delta_{0}}A\to 
A\otimes_{{\mathrm{k}}\Delta_{0}} A_+^{\otimes_{{\mathrm{k}}\Delta_{0}}^{n-2}}\otimes_{{\mathrm{k}}\Delta_{0}}A 
\] 
which is the $A$-bimodule extension of 
\[ 
\phi_{\alpha}(1\otimes v_1\dots v_m\otimes1)\!=\! 
\begin{cases} 
\underbrace{v_1\dots v_{s_1}}_{x_1} 
[ 
\underbrace{\dots v_{s_2}}_{x_2}| 
\dots| 
\underbrace{\dots v_{s_{n-1}}}_{x_{n-1}}] 
\underbrace{\dots v_{m}}_{x_{n}},& 
\text{if $c_j<N$ $\forall j$;} \\ 
0,& 
\text{otherwise.} 
\end{cases} 
\] 
where $s_i=x_1+\dots+x_i$, $i=1\dots n$. 
Thus we obtain for all $m$ a map 
\[ 
\phi:{\bf C}_n(m,N)\to\Hom_{A^e} 
(A\otimes_{\Delta_{0}}{\mathrm{k}}\Delta_{m}\otimes_{\Delta_{0}}A,{\bf Q}_{n-2}). 
\] 
In this context, the comparison morphism ${\bf F}:{\bf P}\to{\bf Q}$ (see \S\ref{subsec:F}) is 
\[ 
{\bf F}_n= 
\begin{cases} 
\phi_{\tilde A_{k(N-1)}^k}: 
{\bf P}_{2k}= 
A\otimes_{\Delta_{0}}{\mathrm{k}}\Delta_{kN}\otimes_{\Delta_{0}}A\longrightarrow{\bf Q}_n, 
&\text{if $n=2k$,}\rule[-10pt]{0pt}{10pt} 
 \\ 
\phi_{\tilde B_{k(N-1)}^k}: 
{\bf P}_{2k+1}= 
A\otimes_{\Delta_{0}}{\mathrm{k}}\Delta_{kN+1}\otimes_{\Delta_{0}}A\longrightarrow{\bf Q}_n, 
&\text{if $n=2k+1$.} 
\end{cases} 
\] 
Now the proof of Theorem \ref{MainComparison} 
will be complete if we show that 
\begin{align*} 
b_{2k+1}\circ\phi_{\tilde B_{k(N-1)}^k}&=\phi_{\tilde A_{k(N-1)}^k}\circ d_{2k+1} 
=\phi_{[0,A_{k(N-1)}^k]_N-[1,A_{k(N-1)}^k]_N}\\ 
\intertext{and} 
b_{2k}\circ\phi_{\tilde A_{k(N-1)}^k}&=\phi_{\tilde B_{(k-1)(N-1)}^{(k-1)}}\circ d_{2k} 
=\phi_{\sum_{j=0}^{N-1}[j,B_{k(N-1)-j}^{k-1}]_{N}}, 
\end{align*} 
 
which is the commutativity of the diagrams 
\[ 
\xymatrix 
{ 
\scriptstyle 
A\otimes_{{}_{\!\Delta_{0}}}\!\! 
A_+^{\otimes^{2k+1}_{\Delta_0}} 
\otimes_{{}_{\!\Delta_{0}}}\!\!\!\!A 
& 
\scriptstyle 
A\otimes_{{}_{\!\Delta_{0}}}\!\! 
A_+^{\otimes^{2k}_{\Delta_0}} 
\otimes_{{}_{\!\Delta_{0}}}\!\!\!\!A 
& 
\scriptstyle 
A\otimes_{{}_{\!\Delta_{0}}}\!\! 
A_+^{\otimes^{2k-1}_{\Delta_0}} 
\otimes_{{}_{\!\Delta_{0}}}\!\!\!\!A 
                                       \\ 
{\bf Q}_{2k+1}\ar@{}|\|[u]\ar[r]^{b_{2k+1}}& 
{\bf Q}_{2k}\ar@{}|\|[u]\ar[r]^{b_{2k}}& 
{\bf Q}_{2k-1}\ar@{}|\|[u] 
                                       \\ 
{\bf P}_{2k+1}\ar[r]_{d_{2k+1}}\ar[u]^{{\bf F}_{2k+1}}& 
{\bf P}_{2k}\ar[r]_{d_{2k}}\ar[u]^{{\bf F}_{2k}}& 
{\bf P}_{2k-1}\ar[u]^{{\bf F}_{2k-1}} 
                                        \\ 
\scriptstyle 
A\otimes_{{}_{\!\Delta_{0}}}\!\!\!\!{\mathrm{k}}\Delta_{k\!N\!+\!1}\otimes_{{}_{\!\Delta_{0}}}\!\!\!\!A 
\ar@{}|\|[u] 
& 
\scriptstyle 
A\otimes_{{}_{\!\Delta_{0}}}\!\!\!\!{\mathrm{k}}\Delta_{k\!N}\otimes_{{}_{\!\Delta_{0}}}\!\!\!\!A 
\ar@{}|\|[u] 
& 
\scriptstyle 
A\otimes_{{}_{\!\Delta_{0}}}\!\!\!\!{\mathrm{k}}\Delta_{(k\!-\!1)\!N\!+\!1}\otimes_{{}_{\!\Delta_{0}}}\!\!\!\!A 
\ar@{}|\|[u] 
} 
\] 
This identities are proved in general in the following proposition (see also Proposition \ref{Prop D}). 
 
\begin{proposition} 
For all $n\ge2$ and $m\ge0$ the following diagram is commutative. 
\[ 
\xymatrix{ 
{\bf C}_{n+1}(m,N) 
\ar[r]^{\textstyle \phi\hspace{2.3cm}}\ar[d]_{D_N} & 
\Hom_{A^e}(A\otimes_{\Delta_{0}}{\mathrm{k}}\Delta_{m}\otimes_{\Delta_{0}}A,{\bf Q}_{n-1}) 
\ar[d]^{\textstyle b\circ\_} \\ 
{\bf C}_n(m,N) 
\ar[r]_{\textstyle \phi\hspace{2.3cm}} & 
\Hom_{A^e}(A\otimes_{\Delta_{0}}{\mathrm{k}}\Delta_{m}\otimes_{\Delta_{0}}A,{\bf Q}_{n-2}) 
} 
\] 
In other words, 
$ 
b\circ\phi_{\alpha}=\phi_{{D_N}(\alpha)} 
$ 
for all 
$\alpha\in{\bf C}_{n+1}(m,N)$. 
\end{proposition} 
 
\begin{proof} 
It is sufficient to prove that 
$ 
b(\phi_{\alpha}(T))= 
\phi_{{D_N}(\alpha)}(T)\rule[-8pt]{0pt}{10pt} 
$ 
for all the monomials $T$ of the form 
$T=1\otimes v_1\dots v_m\otimes1\in 
A\otimes_{\Delta_{0}}{\mathrm{k}}\Delta_{m}\otimes_{\Delta_{0}}A$ 
and for all the compositions 
$\alpha=[x_0,x_1,\dots,x_n]\in {(C^{m}_N)}_{n+1}$ with 
$x_j<N$ for all $j=0\dots n$. 
If we denote by $s_i=x_0+\dots+x_i$, $i=0\dots n$ then 
both sides of the above equality are 
\begin{align*} 
&\underbrace{v_1\dots v_{s_1}}_{x_0+x1} 
[ 
\underbrace{v_{s_1+1}\dots v_{s_2}}_{x_2}| 
\dots| 
\underbrace{v_{s_{n-2}+1}\dots v_{s_{n-1}}}_{x_{n-1}} 
] 
\underbrace{v_{s_{n-1}+1}\dots v_{m}}_{x_{n}} 
                                                   \\ 
&+\sum_{i=1}^{n-1} 
(-1)^i 
\underbrace{v_1\dots v_{s_0}}_{x_0} 
[ 
\underbrace{\dots v_{s_1}}_{x_1}| 
\dots| 
\underbrace{v_{s_{i-1}+1}\dots v_{s_{i+1}}}_{\;\;\;\;x_{i}+x_{i+1}}| 
\dots| 
\underbrace{\dots v_{s_{n-1}}}_{x_{n-1}} 
] 
\underbrace{v_{s_{n-1}+1}\dots v_{m}}_{x_{n}}\rule[-28pt]{0pt}{10pt} 
                                                   \\ 
&+(-1)^{n} 
\underbrace{v_1\dots v_{s_0}}_{x_0} 
[ 
\underbrace{v_{s_0+1}\dots v_{s_1}}_{x_1}| 
\dots| 
\underbrace{v_{s_{n-3}+1}\dots v_{s_{n-2}}}_{x_{n-2}} 
] 
\underbrace{v_{s_{n-2}+1}\dots v_{m}}_{x_{n-1}+x_{n}}, 
\end{align*} 
as it follows from the definition of $\phi$, ${D_N}$ and $b$. 
\end{proof}

\section{The Hochschild cohomology ring} 
\label{sec:cohomology ring} 

We begin this section by recalling some definitions and notation 
following \cite{Ci2}. 
As we said before, all quivers are assumed to be finite and connected. 
 
A path $\gamma$ in a quiver is said to be an {\sl oriented cycle} 
if $o(\gamma)=t(\gamma)$. 
 
Two  paths $\alpha\in\Delta_i$ and $\beta\in\Delta_j$ are {\sl parallel}, 
if $o(\alpha)=o(\beta)$ and $t(\alpha)=t(\beta)$. 
Let 
\[ 
\Delta_i\parallel\Delta_j=\{(\alpha,\beta):\, \alpha\in\Delta_i,\beta\in\Delta_j 
\text{ and $\alpha$ is parallel to $\beta$}\}. 
\] 
 
A pair $(\alpha,\beta)$ of parallel paths is said to {\sl start 
together} if they have the first arrow in common, and they are said 
to {\sl end  together} if they have the last arrow in common. 
 
A vertex is called a {\sl sink} (resp.\ a  {\sl source}) 
if it is not the source (resp.\ end) vertex of any arrow. 
 
Parallel paths that start together and do not end at a sink 
can be pushed forward. 
More precisely, let $(\alpha,\beta)$ be a pair of parallel paths that 
start together. 
Then $\alpha=v\gamma$, $\beta=v\delta$ with $v\in\Delta_1$ and 
$t(\alpha)=t(\beta)$. 
Then any pair $(\tilde \alpha,\tilde \beta)$ satisfying 
$\tilde \alpha=\gamma w$, $\tilde\beta=\delta w$, with $w\in\Delta_1$ and 
$o(w)=t(\alpha)=t(\beta)$ is called a {\sl $+$movement} of $(\alpha,\beta)$. 
In an analogous way are define {\sl $-$movements}. 
 
A pair of parallel paths $(\alpha,\beta)$ is said to be a {\sl $+$extreme} 
({\sl $-$extreme}) if it does not admit any $+$movement ($-$movement). 
Therefore, a pair of parallel paths $(\alpha,\beta)$ is a 
$+$extreme if and only if they end at a sink or 
do not start together (clearly both might occur simultaneously). 
An analogous characterization holds for $-$extremes. 
We shall call a pair of parallel paths $(\alpha,\beta)$ 
just an {\sl extreme} if it is either a $+$extreme or a $-$extreme. 
 
Finally, two pairs $(\alpha,\beta)$ and $(\gamma,\delta)$ in 
$\Delta_i\parallel\Delta_j$ are said to be {\sl equivalent}, 
and denoted by $(\alpha,\beta)\sim(\gamma,\delta)$, 
 if there exist a finite 
sequence of $+$movements and $-$movements carrying 
$(\alpha,\beta)$ to $(\gamma,\delta)$.

\begin{definition} 
An equivalence class in $\Delta_i\parallel\Delta_j$ is called a 
{\sl medal} if all its $+$extremes end at a sink and all its 
$-$extremes start at a source. 
In particular, a class without extremes is a medal. 
\end{definition} 
 
\noindent 
\textbf{Examples.} 
In the first example, let $\alpha=v_1v_2$ and  $\beta=v_1v_2v_3v_4v_1v_2$. 
The class of $(\alpha,\beta)$ is a medal since it does not contain any extreme. 
In fact, any pair of parallel paths $(\alpha,\beta)$ in an oriented cycle 
can be pushed forward and pulled backwards 
and therefore there are no extremes. 
In particular, every class is a medal. 
 
In the second example, let $\alpha=v_1v_2$ and $\beta=v_1v_2v_3v_4v_1v_2$. 
Although $(\alpha,\beta)$ could be pushed forward indefinitely, it is not a medal since 
\[ 
(\alpha,\beta) 
\sim(v_2v_3,v_2v_3v_4v_1v_2v_3) 
\sim(v_3v_2,v_3v_4v_1v_2v_3v_2) 
\sim(v_2v_3,v_4v_1v_2v_3v_2v_3) 
\] 
and the last pair is a $+$extreme that does not end at a sink. 
\begin{center} 
\setlength{\unitlength}{.8cm} 
\begin{picture}(2,3.5)(0,-1) 
  \put(0,2){\circle*{.15}} 
  \put(0.2,1.2){\oval(1.6,1.6)[tr]} 
  \put(0.1,2){\vector(-1,0){0}} 
  \put(1,1){\circle*{.15}} 
  \put(.2,0.8){\oval(1.6,1.6)[br]} 
  \put(1,.9){\vector(0,1){0}} 
  \put(0,0){\circle*{.15}} 
  \put(-.2,0.8){\oval(1.6,1.6)[bl]} 
  \put(-0.1,0){\vector(1,0){0}} 
  \put(-1,1){\circle*{.15}} 
  \put(-.2,1.2){\oval(1.6,1.6)[tl]} 
  \put(-1,1.1){\vector(0,-1){0}} 
  \put(1,1.5){$v_1$} 
  \put(-1.5,1.5){$v_2$} 
  \put(-1.5,0.2){$v_3$} 
  \put(1,0.2){$v_4$} 
  \put(-1,-1){Example 1.} 
\end{picture} 
\hspace{2cm} 
\setlength{\unitlength}{.6cm} 
\begin{picture}(2,2)(0,-1.6) 
  \put(1,1){\circle*{.18}} 
  \put(-0.05,1.3){\oval(1.9,1.5)[t]} 
  \put(-1,1.2){\vector(0,-1){0}} 
  \put(-1,1){\circle*{.18}} 
  \put(-0.05,.7){\oval(1.9,1.5)[b]} 
  \put(0.9,0.8){\vector(0,1){0}} 
  \put(3,1){\circle*{.18}} 
  \put(2.05,1.3){\oval(1.9,1.5)[t]} 
  \put(1.1,1.2){\vector(0,-1){0}} 
  \put(1,1){\circle*{.18}} 
  \put(2.05,.7){\oval(1.9,1.5)[b]} 
  \put(3,0.8){\vector(0,1){0}} 
  \put(1.8,2.3){$v_1$} 
  \put(-0.4,2.3){$v_2$} 
  \put(-0.4,-0.6){$v_3$} 
  \put(1.8,-0.6){$v_4$} 
  \put(-.8,-1.6){Example 2.} 
\end{picture} 
\end{center} 
 
\ 
 
\subsection{The Hochschild cohomology groups} 

The zero cohomology group $H^0(A,A)$ is the center of $A$ as for any algebra. 
We now describe the cohomology groups $H^n(A,A)$, for $n>0$, 
following \cite{Lo}. 
Being 
\[ 
{\bf P}_n= 
\begin{cases} 
A\otimes_{\Delta_{0}}{\mathrm{k}}\Delta_{kN}\otimes_{\Delta_{0}}A, &\text{if $n=2k$;} \\ 
A\otimes_{\Delta_{0}}{\mathrm{k}}\Delta_{kN+1}\otimes_{\Delta_{0}}A, &\text{if $n=2k+1$;} 
\end{cases} 
\] 
we have 
\begin{equation}\label{HomP1} 
\Hom_{A^e}({\bf P}_n,A)\simeq 
\begin{cases} 
\Hom_{\Delta_{0}^e}({\mathrm{k}}\Delta_{kN},A), &\text{if $n=2k$;} \\ 
\Hom_{\Delta_{0}^e}({\mathrm{k}}\Delta_{kN+1},A), &\text{if $n=2k+1$.} 
\end{cases} 
\end{equation} 
Let ${\bf P}^*$ be the bigraded vector space 
\[ 
{\bf P}^*=\bigoplus_{n\ge 0}\bigoplus_{i=0}^{N-1}{\bf P}^*_{n,i} 
\] 
where 
\[ 
{\bf P}^*_{2k,i}={\mathrm{k}}\Delta_i\parallel\Delta_{kN} 
\qquad\text{and}\qquad 
{\bf P}^*_{2k+1,i}={\mathrm{k}}\Delta_i\parallel\Delta_{kN+1} 
\] 
Since it is clear that 
$\Hom_{\Delta_{0}^e}({\mathrm{k}}\Delta_{m},A) 
\simeq\bigoplus_{j=0}^{N-1}{\mathrm{k}}\Delta_j\parallel\Delta_m$ for all $m\ge0$, 
it follows that 
\begin{equation}\label{HomP2} 
\Hom_{A^e}({\bf P}_n,A)\simeq \bigoplus_{i=0}^{N-1}{\bf P}^*_{n,i}. 
\end{equation} 
The following theorem is proved in \cite[\S3]{Lo} and 
it describes the cohomology of truncated quiver algebras. 
We note the word {\sl $j$-extreme} is used instead of {\sl medal} in Locateli's paper. 
 
\begin{theorem}\label{Thm:Lo} 
Let $\Delta$ be a quiver. 
Then the complex $\Hom_{A^e}({\bf P}_n,A)$ 
has the following decomposition into subcomplexes 
\begin{center} 
\setlength{\unitlength}{0.94cm} 
\begin{picture}(10,7)(1.2,-7.4) 
  \put(-.4,-1){$\Hom_{\Delta_0^e}({\mathrm{k}}\Delta_{kN},A)$} 
  \put(2.8,-.9){\vector(1,0){1.1}} 
  \put(3.1,-0.75){$\scriptstyle d_{2k}$} 
  \put(4.2,-1){$\Hom_{\Delta_0^e}( {\mathrm{k}}\Delta_{kN+1},A)$} 
  \put(7.8,-.9){\vector(1,0){1.1}} 
  \put(8,-0.75){$\scriptstyle d_{2k+1}$} 
  \put(9.1,-1){$\Hom_{\Delta_0^e}( {\mathrm{k}}\Delta_{(k+1)N},A)$} 
\put(-.5,-1.3){\line(1,0){13.1}} 
  \put(0,-2){${\mathrm{k}}\;\Delta_0\parallel\Delta_{ kN}$} 
  \put(4.6,-2){${\mathrm{k}}\;\Delta_0\parallel\Delta_{ kN+1}$} 
  \put(9.5,-2){${\mathrm{k}}\;\Delta_0\parallel\Delta_{ (k+1)N}$} 
  \put(0,-3){${\mathrm{k}}\;\Delta_1\parallel\Delta_{ kN}$} 
  \put(4.6,-3){${\mathrm{k}}\;\Delta_1\parallel\Delta_{ kN+1}$} 
  \put(9.5,-3){${\mathrm{k}}\;\Delta_1\parallel\Delta_{ (k+1)N}$} 
  \put(0,-4){${\mathrm{k}}\;\Delta_2\parallel\Delta_{ kN}$} 
  \put(4.6,-4){${\mathrm{k}}\;\Delta_2\parallel\Delta_{ kN+1}$} 
  \put(9.5,-4){${\mathrm{k}}\;\Delta_2\parallel\Delta_{ (k+1)N}$} 
  \put(1,-4.8){$\cdot$} 
  \put(6,-4.8){$\cdot$} 
  \put(10.5,-4.8){$\cdot$} 
  \put(1,-5){$\cdot$} 
  \put(6,-5){$\cdot$} 
  \put(10.5,-5){$\cdot$} 
  \put(1,-5.2){$\cdot$} 
  \put(6,-5.2){$\cdot$} 
  \put(10.5,-5.2){$\cdot$} 
  \put(0,-6){${\mathrm{k}}\;\Delta_{N-2}\parallel\Delta_{ kN}$} 
  \put(4.6,-6){${\mathrm{k}}\;\Delta_{N-2}\parallel\Delta_{ kN+1}$} 
  \put(9.5,-6){${\mathrm{k}}\;\Delta_{N-2}\parallel\Delta_{ (k+1)N}$} 
  \put(0,-7){${\mathrm{k}}\;\Delta_{N-1}\parallel\Delta_{ kN}$} 
  \put(4.6,-7){${\mathrm{k}}\;\Delta_{N-1}\parallel\Delta_{ kN+1}$} 
  \put(9.5,-7){${\mathrm{k}}\;\Delta_{N-1}\parallel\Delta_{ (k+1)N}$} 
  \put(2.4,-2.2){\vector(3,-1){1.8}}\put(3.2,-2.3){$\scriptstyle D_0^{2k}$} 
  \put(2.4,-3.2){\vector(3,-1){1.8}}\put(3.2,-3.3){$\scriptstyle D_1^{2k}$} 
  \put(3.3,-4.8){$\cdot$} 
  \put(3.3,-5.0){$\cdot$} 
  \put(3.3,-5.2){$\cdot$} 
  \put(2.4,-6.2){\vector(3,-1){1.8}}\put(3.2,-6.3){$\scriptstyle D_{N-2}^{2k}$} 
  \put(8.5,-5){$\scriptstyle D_{0}^{2k+1}$} 
  \put(7.6,-1.9){\line(1,0){.5}} 
  \put(8.1,-2.15){\oval(0.5,0.5)[tr]} 
  \put(8.35,-2.15){\line(0,-1){4.5}} 
  \put(8.6,-6.65){\oval(0.5,0.5)[bl]} 
  \put(8.6,-6.9){\vector(1,0){.5}} 
%
\end{picture} 
\end{center} 
where 
\begin{align*} 
D_j^{2k}(\alpha,\pi)&=\sum_{a\in\Delta_1}(a\alpha,a\pi)-\sum_{b\in\Delta_1}(\alpha b,\pi b), 
\;k\ge0\text{ and $j=0,\dots,N-2$}, \\ 
D_0^{2k+1}(v,\pi)&=\sum_{ab\in\Delta_{N-1}}(avb,a\pi b), 
\;k\ge0. 
\end{align*} 
The following holds for the differentials. 
\begin{enumerate}[1.] 
\item $D_0^{2k+1}$ is injective for all $k\ge0$, 
\item $D_0^{2k}$ is injective for all $k>0$ unless $\Delta$ is a cycle, 
\item $\dim\ker D_j^{2k}$ is equal to the number of medals in 
$\Delta_{j}\parallel\Delta_{ kN}$ for all $j=1,\dots,N-2$ and for all $k>0$. 
More precisely, if for each medal $M$ one considers 
\[ 
\bar M=\sum_{(\alpha,\pi)\in M} (\alpha,\pi)\in\mathrm{k}\,\Delta_{j}\parallel\Delta_{ kN}, 
\] 
then the set $\{\bar M: M\text{ is a medal in $\Delta_{j}\parallel\Delta_{ kN}$}\}$ 
is a basis of $\ker D_j^{2k}$. 
\end{enumerate} 
\end{theorem} 
\begin{definition}\label{Def:medal class} 
The cohomology class $\bar M$ corresponding to the medal $M$ 
will be called the \textit{medal cohomology class} associated to $M$. 
\end{definition} 
 
\begin{remark}\label{rmk:bigraded} 
The cohomology group $H^*(A,A)$ inherits the bigrading of the complex 
${\bf P}^*$. We now have 
 \begin{enumerate}[1.] 
  \item $H^n_0=0$, for all $n\ge 1$. 
  \item $H^{2k}_i$ is formed entirely of medals for all $1\le i\le N-2$. 
  \item $H^{2k}_{N-1}$ is a cokernel. 
  \item $H^{2k+1}_i$ is a cokernel for all $1\le i\le N-1$. 
 \end{enumerate} 
\end{remark}

\subsection{The Yoneda product} 
The Hochschild cohomology groups $H^*(A,A)$ of $A$ are, by definition, 
$\text{Ext}_{A^e}^*(A,A)$ and therefore they have a ring structure 
given by the multiplication induced by the Yoneda product. 
 
It is well known that the Yoneda product of $H^*(A,A)$ coincides with the cup product 
defined on the cohomology of $\text{Hom}_{A^e}({\bf Q},A)$. 
More precisely, the cup product is originally defined in terms of the 
standard $A^e$-projective bar resolution $A\otimes A^{\otimes *}\otimes A$ 
of $A$ as follows. 
Given two cochains, 
\begin{align*} 
f\in\Hom_{A^e}(A\otimes A^{\otimes m}\otimes A,A)&\simeq\Hom_{\mathrm{k}}(A^{\otimes m},A) \\ 
g\in\Hom_{A^e}(A\otimes A^{\otimes n}\otimes A,A)&\simeq\Hom_{\mathrm{k}}(A^{\otimes n},A) 
\end{align*} 
the cup product of $f$ and $g$ is the cochain 
$ 
f\cup g\in\Hom_{\mathrm{k}}(A^{\otimes (m+n)},A) 
$ 
defined by 
\[ f\cup g\;(\alpha_1\otimes\dots\otimes\alpha_{m+n})= 
f(\alpha_1\otimes\dots\otimes\alpha_{m})\,g(\alpha_{m+1}\otimes\dots\otimes\alpha_{m+n}). 
   \] 
 
The analogous definition works for the resolution ${\bf Q}$. 
For any ${\mathrm{k}}\Delta_0$-bimodule $M$, 
let $\Hom_{\Delta_0^e}(M,A)$ be the group of 
homomorphisms of ${\mathrm{k}}\Delta_0$-bimodules. 
Given two cochains, 
\begin{align*} 
f\in\Hom_{A^e}({\bf Q}_m,A)&\simeq\Hom_{\Delta_0^e}(A_+^{\otimes^m_{\Delta_0}},A) \\ 
g\in\Hom_{A^e}({\bf Q}_n,A)&\simeq\Hom_{\Delta_0^e}(A_+^{\otimes^n_{\Delta_0}},A) 
\end{align*} 
the cup product of $f$ and $g$ is the cochain 
$ 
f\cup g\in\Hom_{\Delta_0^e}(A_+^{\otimes^{m+n}_{\Delta_0}},A) 
$ 
defined by 
\begin{equation}\label{eq:cup} 
 f\cup g\;([\alpha_1|\dots|\alpha_{m+n}])= 
f([\alpha_1|\dots|\alpha_{m}])\,g([\alpha_{m+1}|\dots|\alpha_{m+n}]). 
\end{equation} 
 
We shall now use the comparison morphisms to 
describe the Yoneda product in the minimal resolution ${\bf P}$ of $A$.

\begin{proposition} 
Let $f\in\Hom_{A^e}({\bf P}_m,A)$ and $g\in\Hom_{A^e}({\bf P}_n,A)$. 
Then, in terms of the identification (\ref{HomP1}), we have 
\begin{enumerate}[$\bullet$] 
\item 
if $m=2h$ and $n=2k$ then 
\[ 
f\cup g\,(v_1\dots v_{hN}w_1\dots w_{kN})= 
f(v_1\dots v_{hN})g(w_1\dots w_{kN}), 
\] 
\item 
if $m=2h$ and $n=2k+1$ then 
\[ 
f\cup g\,(v_1\dots v_{hN}w_1\dots w_{kN})= 
f(v_1\dots v_{hN})g(w_1\dots w_{kN+1}), 
\] 
\item 
if $m=2h+1$ and $n=2k+1$ then 
\begin{multline*} 
f\cup g\,(u_1\dots u_{(h+k+1)N})= \\ 
\sum_{0<i<j<N} 
u_1\dots f(u_i\dots u_{i+hN})u_{i+hN+1}\dots 
g(u_{j+hN}\dots u_{j+(h+k)N})\dots u_{(h+k+1)N} 
\end{multline*} 
\end{enumerate} 
\end{proposition} 
 
\begin{proof} 
By definition $f\cup g=((f\circ\mathbf{G})\cup (g\circ\mathbf{G}))\circ\mathbf{F}$. 
 
Assume that $m=2h$ and $n=2k$ and let 
$u_1\dots u_{lN}=v_1\dots v_{hN}w_1\dots w_{kN}$, $l=h+k$. 
Then 
\begin{multline*} 
{\bf F}(u_1\dots u_{lN})= 
{\bf F}(1\otimes u_1\dots u_{lN}\otimes1)= \\ 
\sum 
1 
[ 
\underbrace{u_1\dots u_{x_1}}_{x_1} 
| 
\underbrace{u_{1+x_1}}_{1} 
| 
\underbrace{\dots u_{1+x_1+x_2}}_{x_2} 
| 
\underbrace{u_{2+x_1+x_2}}_{1} 
|\dots\dots| 
\underbrace{u_{l+\sum x_j}}_{1} 
] 
\underbrace{\dots u_{lN}}_{\!\!\!lN-l-\sum x_j}, 
\end{multline*} 
where the sum is over all $l$-tuples $(x_1,\dots,x_l)\in\Z^l$ 
with $1\le x_i<N$. 
Being 
$f\circ\mathbf{G}\in\Hom_{\Delta_0^e}(A_+^{\otimes^{2h}_{\Delta_0}},A)$ and 
$g\circ\mathbf{G}\in\Hom_{\Delta_0^e}(A_+^{\otimes^{2k}_{\Delta_0}},A)$, 
and identifying $[\alpha_1|\dots|\alpha_{m}]$ with $1[\alpha_1|\dots|\alpha_{m}]1$, 
we have (see \eqref{eq:cup}) 
\begin{multline*} 
((f\circ\mathbf{G})\cup (g\circ\mathbf{G}))\circ\mathbf{F}(u_1\dots u_{lN})= \\ 
\sum 
(f\circ\mathbf{G})\big( 
1 
[ 
\underbrace{u_1\dots u_{x_1}}_{x_1} 
| 
\underbrace{u_{1+x_1}}_{1} 
| 
\dots\dots 
| 
\underbrace{\dots u_{h-1+\sum_{j=1}^h x_j}}_{x_h} 
| 
\underbrace{u_{h+\sum_{j=1}^h x_j}}_{1}]1 
\big) \\ 
\times 
(g\circ\mathbf{G})\big( 
1 
[ 
\underbrace{\dots u_{h+\sum_{j=1}^{h+1} x_j}}_{x_{h+1}} 
| 
\underbrace{u_{h+1+\sum_{j=1}^{h+1} x_j}}_{1} 
|\dots\dots| 
\underbrace{u_{l+\sum x_j}}_{1} 
]1 
\big) 
\underbrace{\dots u_{lN}}_{\!\!\!lN-l-\sum x_j}. 
\end{multline*} 
By the definition of $\mathbf{G}$, the only non zero term of this sum is 
the corresponding to $(x_1,\dots,x_l)=(N-1,N-1,\dots,N-1)$. 
This yields 
\begin{align*} 
f\cup g\,(u_1\dots u_{lN}) 
&=f(u_1\dots u_{hN})g(u_{hN+1}\dots u_{lN})\\ 
&=f(v_1\dots v_{hN})g(w_1\dots w_{kN}). 
\end{align*} 
The proof is analogous if $m=2h$ and $n=2k+1$. 
Finally assume that $m=2h+1$ and $n=2k+1$. 
Let $l=h+k+1$. 
Then 
\begin{multline*} 
((f\circ\mathbf{G})\cup (g\circ\mathbf{G}))\circ\mathbf{F}(u_1\dots u_{lN})= \\ 
\sum 
(f\circ\mathbf{G})( 
1 
[ 
\underbrace{u_1\dots u_{x_1}}_{x_1} 
| 
\underbrace{u_{1+x_1}}_{1} 
| 
\dots\dots 
| 
\underbrace{u_{h+\sum_{j=1}^h x_j}}_{1} 
| 
\underbrace{\dots u_{h+\sum_{j=1}^{h+1} x_j}}_{x_{h+1}}]1)\\ 
\times 
(g\circ\mathbf{G})\big( 
1 
[ 
\underbrace{u_{h+1+\sum_{j=1}^{h+1} x_j}}_{1} 
| 
\underbrace{\dots u_{h+1+\sum_{j=1}^{h+2} x_j}}_{x_{h+2}} 
| 
\dots\dots| 
\underbrace{u_{l+\sum x_j}}_{1} 
]1 
\big) 
\underbrace{\dots u_{lN}}_{\!\!\!lN-l-\sum x_j}. 
\end{multline*} 
Again, by the definition of $\mathbf{G}$, the only non zero terms of this sum are 
the corresponding to 
$(x_1,\dots,x_l)=(x_1,N-1,N-1,\dots,N-1)$ for all $x_1=1,\dots,N-1$. 
Therefore 
\begin{multline*} 
f\cup g\,(u_1\dots u_{lN})= \\ 
\sum_{x_1=1}^{N-1} 
\sum_{j=1}^{x_1} 
u_1\dots f(\underbrace{u_j\dots }_{hN+1})u_{j+hN+1}\dots 
g(\underbrace{u_{x_1+1+hN}\dots }_{kN+1})\dots u_{(h+k+1)N} 
\end{multline*} 
as claimed. 
\end{proof} 
Recall that 
\[ 
{\bf P}^*_{2k,i}={\mathrm{k}}\Delta_i\parallel\Delta_{kN} 
\qquad\text{and}\qquad 
{\bf P}^*_{2k+1,i}={\mathrm{k}}\Delta_i\parallel\Delta_{kN+1} 
\] 
and that we have an isomorphism 
$ 
\Hom_{A^e}({\bf P}_n,A)\simeq \bigoplus_{i=0}^{N-1}{\bf P}^*_{n,i} 
$ 
(see \eqref{HomP2}). 
 
\begin{theorem}\label{Thm:YonedaResolution} 
Let $A$ be a truncated quiver algebra. Then, in terms of the 
above isomorphism we have that 
 
\vspace{2mm} 
 
\begin{enumerate}[$\bullet$] 
\item 
if $(\alpha,\pi)\in{\bf P}^*_{2h,i}$ and 
$(\beta,\tau)\in{\bf P}^*_{2k,j}$, 
then 
\[ 
(\alpha,\pi)\cup(\beta,\tau)=(\alpha\beta,\pi\tau)\in 
{\bf P}^*_{2(h+k),i+j}. 
\] 
 
\vspace{2mm} 
 
\item 
if $(\alpha,\pi)\in{\bf P}^*_{2h,i}$ and 
$(\beta,\tau)\in{\bf P}^*_{2k+1,j}$, 
then 
\[ 
(\alpha,\pi)\cup(\beta,\tau)=(\alpha\beta,\pi\tau)\in 
{\bf P}^*_{2(h+k)+1,i+j}. 
\] 
 
\vspace{2mm} 
 
\item 
if $(\alpha,\pi)\in{\bf P}^*_{2h+1,i}$ and 
$(\beta,\tau)\in\in{\bf P}^*_{2k+1,j}$, 
then 
\[ 
(\alpha,\pi)\cup(\beta,\tau)=\sum_\mu (\gamma_\mu,\mu)\in 
{\bf P}^*_{2(h+k)+2,N-2+i+j}, 
\] 
where the sum is over all paths $\mu$ containing $\pi$ and $\tau$ as a subpath, 
and $\gamma_\mu$ is the result of substituting $\pi$ and $\tau$ 
by $\alpha$ and $\beta$ respectively in $\mu$. 
In particular, 
$(\alpha,\pi)\cup(\beta,\tau)=0$, if $i+j>1$. 
\end{enumerate} 
\end{theorem}

\begin{remark} 
In fact the third line of the previous theorem shows that the Yoneda product 
on the complex ${\bf P}^*$ is not compatible with the bigrading. 
\end{remark}

Let us define the following product in ${\bf P}^*$. 
If $(\alpha,\pi)\in{\bf P}^*_{n_1}$ and $(\beta,\tau)\in{\bf P}^*_{n_2}$ 
let 
\[ 
(\alpha,\pi)\vee(\beta,\tau)= 
\begin{cases} 
(\alpha\beta,\pi\tau),&\text{ if $n_1$ or $n_2$ are even;}\\ 
0,&\text{ otherwise}. 
\end{cases} 
\] 
The following theorem extends the results of Sections 3 and 4 of \cite{BLM}, 
proved for truncated cycle algebras, 
to any truncated quiver algebra. 
 
\begin{theorem}\label{Thm:YonedaCohomology} 
Let $A$ be a truncated quiver algebra. 
Then the product $\vee$ in ${\bf P}^*$ induces the Yoneda product in the 
Hochschild cohomology group $H^*(A,A)$, and thus it is a bigraded 
commutative ring (see Remark \ref{rmk:bigraded}). 
In particular the Yoneda product of two odd-degree cohomology classes is zero. 
\end{theorem} 
\begin{proof} 
If $n_1$ or $n_2$ are even, then the result is straightforward 
from the previous theorem. 
If $n_1=hN+1$ and $n_2=kN+1$ are odd numbers then 
$(\alpha,\pi)\in{\mathrm{k}}\Delta_{i}\parallel\Delta_{ hN+1}$ and 
$(\beta,\tau)\in{\mathrm{k}}\Delta_{j}\parallel\Delta_{ kN+1}$ 
with $1\le i,j\le N-1$ (see Theorem \ref{Thm:Lo}). 
Thus $i+j\ge2$ and, according to the previous theorem, we obtain $(\alpha,\pi)\cup(\beta,\tau)=0$. 
\end{proof} 
 
\begin{corollary}\label{coro:YonedaMedals} 
Let $A$ be a truncated quiver algebra. 
If the Yoneda product of two cohomology classes of positive degree is not zero then 
at least one of them is a medal cohomology class. 
\end{corollary} 
\begin{proof} 
Let $f\in\kk\Delta_{i}\parallel\Delta_{m_1}$ and $g\in\kk\Delta_{j}\parallel\Delta_{m_2}$ be 
representatives of cohomology classes $\bar f$ and $\bar g$ 
of cohomological degrees $n_1>0$ and $n_2>0$ respectively 
and assume that $\bar f\cup \bar g\ne0$. 
Combining Theorem \ref{Thm:Lo} and 
Theorem \ref{Thm:YonedaCohomology} we obtain that either $n_1$ or $n_2$ is even 
and that $i+j\le N-1$, $i>0$, $j>0$. 
Thus $i,j\le N-2$ and assuming that $n_1$ is even 
it follows that $\bar f$ is a medal cohomology class. 
\end{proof}

\begin{corollary} 
If $f_1,\dots,f_N$ are cohomology classes of positive degree, then $f_1\dots f_N=0$. 
In particular $H^*(A,A)/\mathcal{N}\simeq \kk$, where $\mathcal{N}$ is the ideal 
generated by homogeneous nilpotent elements. 
\end{corollary} 
\begin{proof} 
The result follows directly from Theorems \ref{Thm:Lo} and \ref{Thm:YonedaCohomology}. 
\end{proof}

\section{Applications} 
\label{sec:examples} 

Even though it does not appear in the literature as a conjecture, many people believed 
that for TQA's $A$ the cohomology ring should be trivial, except for cycle algebras, 
meaning more precisely that the product in the subring $\bigoplus_{n\ge1}H^n(A,A)$ should be trivial. 
 
After Corollary \ref{coro:YonedaMedals}, the understanding of medals in $\Delta$ is highly 
relevant to determine whether the product in cohomology is trivial or not. 
 
A full classification of TQA's with trivial cohomology ring requires a deeper 
understanding of the spaces of paths, parallel paths and medals of quivers. 
We intend to carry out this classification in a future work. 
 
We present here two large classes of quiver whose associated TQA's have trivial cohomology rings. 
Namely the class of quivers with no oriented cycles and the class of quiver 
with neither sinks nor sources. 
On the other hand we present an interesting example of a small quiver yielding TQA's 
with a non trivial cohomology ring.

\subsection{Quivers with no oriented cycles} 
 
\begin{theorem}\label{thm:no-oriented-cycles} 
Let $\Delta$ be a quiver without any cycle. 
Then the Yoneda product in $\bigoplus_{n\ge1}H^n(A,A)$ is zero. 
\end{theorem} 
 
\begin{proof} 
The set of vertices $\Delta_0$ of a quiver without any cycle is a partial ordered set: 
$v_1\preceq v_2$ if and only if there 
exists a path $\alpha$ such that $o(\alpha)=v_1$ and $t(\alpha)=v_2$. 
 
We first prove that 
for any non zero medal cohomology class $M$ of positive cohomological degree, 
there exist vertices $v_1\prec v_2$ such that 
$o(\beta)\preceq v_1$ and $v_2\preceq t(\beta)$ 
for any pair $(\beta,\tau)\in M$. 
Let $(\alpha,\pi)$ be any pair in $M$. 
We know that $1\le|\alpha|\le N-2<N\le|\pi|$ (cf.\ Theorem \ref{Thm:Lo}). 
Assume that $\alpha=a_1\dots a_{|\alpha|}$ and $\pi=p_1\dots p_{|\pi|}$. 
Let $l\ge0$ be the largest integer for which $a_i=p_i$ for $i=1,\dots, l$. 
If $l=|\alpha|$ then $p_{l+1}\dots p_{|\pi|}$ would be an oriented cycle, 
which is imposible, and thus $l<|\alpha|$. 
Similarly, if $r\ge0$ is the smallest integer for which 
$a_i=p_{|\pi|-|\alpha|+i}$ for $i=r,\dots, |\alpha|$ 
then $r>1$. Clearly $l<r$ and in fact $l< r-1$ 
for if $l=r-1$ then $p_{l+1}\dots p_{|\pi|-|\alpha|+r-1}$ would be an oriented cycle. 
Therefore, if $v_1=t(a_l)$ and $v_2=o(a_r)$, then $v_1\prec v_2$. 
Since $a_{l+1}\ne p_{l+1}$ and $a_{r-1}\ne p_{|\pi|-|\alpha|+r-1}$ it follows 
that  for any pair $(\beta,\tau)\in M$, $\beta$ must contain the path 
$a_{l+1}\dots a_{r-1}$. Therefore $o(\beta)\preceq v_1$ and $v_2\preceq t(\beta)$. 
 
We now prove that the Yoneda product in $\oplus_{n\ge1}H^n(A,A)$ is zero. 
Assume, on the contrary, that there are 
two cohomology classes of positive cohomological degree $M_1$ and $M_2$ 
such that $M_1\cup M_2\ne0$. 
By Corollary \ref{coro:YonedaMedals} we may assume that $M_1$ is a 
medal cohomology class. Let $v_1\prec v_2$ as above. 
We now must consider two possibilities: 
\begin{enumerate}[1.] 
  \item $M_2$ is also a medal cohomology class. Let $w_1\prec w_2$ as above. 
Since $M_1\cup M_2\ne0$ 
 there exist $(\alpha,\pi)\in M_1$ and $(\beta,\tau)\in M_2$ 
such that $v_2\preceq t(\alpha)=o(\beta)\preceq w_1$. 
Additionally, since the cup product is commutative, 
there exist $(\alpha',\pi')\in M_1$ and $(\beta',\tau')\in M_2$ 
such that $w_2\preceq t(\beta')=o(\alpha')\preceq v_1$. 
This is a contradiction since $v_1\prec v_2$ and $w_1\prec w_2$. 
\item $M_2$ ia an odd cohomology class. 
Since $M_2\cup M_1=M_1\cup M_2\ne0$ 
 there exist a representative $(\beta,\tau)$ of $M_2$ 
and pairs $(\alpha,\pi),(\alpha',\pi')\in M_1$ such that 
$v_2\preceq t(\alpha)=o(\beta)\preceq t(\beta)\preceq o(\alpha)\prec v_1$. 
which is a contradiction. 
\end{enumerate} 
This completes the proof. 
\end{proof} 
 
\begin{remark} 
Quivers in this class might have lots of medals. 
\end{remark}

\subsection{Quivers with neither sinks nor sources} 
 
Truncated tensor algebras are particular cases of truncated quiver algebras 
associated to quivers without sinks and sources. 
Indeed, let $V$ be a finite dimensional $\kk$-vector space and let 
$A_N=T(V)/(V^{\otimes N})$ the $N$-truncated tensor algebra. 
Then $A_N$ is a truncated quiver algebra 
corresponding to the quiver 
 
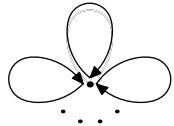
\begin{center} 
   \includegraphics{quiver.1} 
\end{center} 
with $\dim_{\kk}V$ loops. 
If $\dim_{\kk}V\ge2$ then it is not difficult to see that this quiver has no medals in 
$\Delta_j\parallel\Delta_{kN}$ for $k>0$. 
According to Corollary \ref{coro:YonedaMedals} the 
Yoneda product in positive cohomology degrees is zero. 
In Theorem \ref{Thm:TrivialCohomology} below we extend this result to 
all quivers without sinks and sources that are not an oriented cycle. 
We point out that if $\dim_{\kk}V=1$ then the $N$-truncated tensor algebra 
$A_N$ is a truncated cycle algebra 
and its ring structure is described in Subsection \ref{subsec:Constructing}.

\begin{lemma}\label{lemma:quiver1} 
Let $\Delta$ be a quiver. If $(\alpha,\beta)$ is a pair of 
parallel paths such that $|\alpha|<|\beta|$ and its class 
contains no $+$extremes, then there exists an oriented cycle $\gamma$ 
such that $\beta=\alpha\gamma$. Similarly, if $(\alpha,\beta)$ is 
a pair of parallel paths such that $|\alpha|<|\beta|$ and its 
class contains no $-$extremes, then there exists an oriented cycle 
$\gamma$ such that $\beta=\gamma\alpha$. 
\end{lemma} 
\begin{proof} 
Let 
$(\alpha,\beta)$ be a pair of parallel paths such that $|\alpha|<|\beta|$. 
Then 
$\alpha=a_1\dots a_{|\alpha|}$ and 
$\beta=b_1\dots b_{|\alpha|}\dots b_{|\beta|}$. 
If the class of $(\alpha,\beta)$ contains no $+$extremes 
then, in particular, $(\alpha,\beta)$ can be pushed forward $|\alpha|$ times. 
Therefore $b_j=a_j$ for all $j=1,\dots,|\alpha|$. 
Let $\gamma=b_{|\alpha|+1}\dots b_{|\beta|}$. 
Then $o(\gamma)=t(\alpha)=t(\beta)=t(\gamma)$ and thus $\gamma$ is an oriented cycle. 
 
The proof of the second statement is analogous. 
\end{proof} 
 
\begin{lemma}\label{lemma:quiver2} 
Let $\Delta$ be a quiver. 
Then $\Delta$ is an oriented cycle if and only if 
there exist an oriented cycle $\gamma$ and a 
pair of parallel paths $(\alpha,\beta)$, with $\alpha$ and $\beta$ subpaths of $\gamma$ 
and $|\alpha|<|\beta|$, 
such that the class of $(\alpha,\beta)$ does not have any extreme. 
\end{lemma} 
\begin{proof} 
If $\Delta$ is an oriented cycle then it is clear that there are no extremes at all. 
 
Conversely, assume that $\Delta$ is not an oriented cycle. 
We shall prove that given an oriented cycle $\gamma$ 
and a pair of parallel paths $(\alpha,\beta)$, 
with $\alpha$ and $\beta$ subpaths of $\gamma$ 
and $|\alpha|<|\beta|$, 
the class of $(\alpha,\beta)$ has an extreme. 
 
Let $\gamma$ be the oriented cycle. 
Since $\Delta$ is not an oriented cycle there must exist a vertex $p$ in 
$\gamma$ and two different arrows $u,v\in\Delta_1$ (not necessarily in $\gamma$) 
such that either $t(u)=t(v)=p$ or $o(u)=o(v)=p$.

Now let $(\alpha,\beta)$ be a pair of parallel paths such that 
$\alpha$ and $\beta$ are subpaths of $\gamma$ and $|\alpha|<|\beta|$. 
 
In the case $o(u)=o(v)=p$ we shall see that $(\alpha,\beta)$ can be pushed 
forward until we obtain a $+$extreme. 
We first push  $(\alpha,\beta)$ forward 
in order to reach a pair $(\alpha',\beta')$ with $t(\alpha')=t(\beta')=p$. 
Since $\gamma$ is an oriented cycle this is possible, unless we reach a $+$extreme before. 
If $(\alpha',\beta')$ is a $+$extreme we are done. 
Otherwise $\alpha'$ and $\beta'$ must start together and hence 
$\alpha'=a'_1\dots a'_{|\alpha|}$ and $\beta'=b'_1\dots b'_{|\beta|}$ 
with $a'_1=b'_1$. 
Since $|\alpha|+1\le|\beta|$ and  $t(\alpha')=p$ we can assume that $u\ne b'_{|\alpha|+1}$. 
Next, we push $(\alpha',\beta')$ forward obtaining 
\[ 
(\alpha',\beta')\sim 
(a'_2\dots a'_{|\alpha|}u,b'_2\dots b'_{|\alpha|} b'_{|\alpha|+1}\dots b'_{|\beta|}u). 
\] 
Now, if we keep pushing forward, 
since $u\ne b'_{|\alpha|+1}$ it is clear that in at most $|\alpha|-1$ times 
we shall reach a $+$extreme. 
 
In the case $t(u)=t(v)=p$ an analogous argument shows how to 
pull $(\alpha,\beta)$ backwards until a $-$extreme is reached. 
\end{proof}

\begin{theorem}\label{thm:extremes} 
Let $\Delta$ be a finite connected quiver that is not an oriented cycle, 
and let $(\alpha,\beta)$ be a pair of parallel paths such that $|\alpha|<|\beta|$. 
Then the class of $(\alpha,\beta)$ has an extreme. 
In particular, if in addition $\Delta$ has 
neither sinks nor sources then $\Delta_i\parallel\Delta_j$ does not have any medal for all 
$i\ne j$. 
\end{theorem} 
 
\begin{proof} 
Let $(\alpha,\beta)$ be a pair of parallel paths such that $|\alpha|<|\beta|$ 
and assume its class of $(\alpha,\beta)$ contains no $+$extremes. 
Then, by Lemma \ref{lemma:quiver1} there exists an oriented cycle $\gamma=v_1\dots v_k$, $k\ge1$, 
such that $\beta=\alpha\gamma$. 
Let $v_{k+j}=v_j$ for all $j>0$. 
Thus, by pushing forward, 
\[ 
(\alpha,\beta) 
=(\alpha,\alpha v_1\dots v_k) 
\sim(v_{k+1}\dots v_{k+|\alpha|},v_1\dots v_kv_{k+1}\dots v_{k+|\alpha|}) 
=(\tilde\alpha,\tilde\beta). 
\] 
Now $(\tilde\alpha,\tilde\beta)$ is a pair of parallel paths contained in 
the oriented cycle $\gamma$. 
Since $|\tilde\alpha|=|\alpha|<|\beta|=|\tilde\beta|$ and $\Delta$ 
is not an oriented cycle, Lemma \ref{lemma:quiver2} 
implies that there exist an extreme pair in the class of $(\tilde\alpha,\tilde\beta)$. 
 
In particular, if $\Delta$ has no sinks and no sources and 
$(\alpha,\beta)\in\Delta_i\parallel\Delta_j$ with  $i\ne j$ then $(\alpha,\beta)$ 
is equivalent to either a $+$extreme not ending at a sink or a 
$-$extreme not starting at a source. 
Thus the class of $(\alpha,\beta)$ is not a medal. 
\end{proof}

\begin{theorem}\label{Thm:TrivialCohomology} 
Let $\Delta$ be a quiver that is not an  oriented cycle and with 
neither sinks nor sources. 
Then the Yoneda product in $\bigoplus_{n\ge1}H^n(A,A)$ is zero. 
\end{theorem} 
 
\begin{proof} 
From Theorem \ref{thm:extremes} we know that 
$\Delta_{i}\parallel\Delta_{j}$ contains no medals when $i\ne j$ 
and therefore Corollary \ref{coro:YonedaMedals} implies 
that the Yoneda product in $\bigoplus_{n\ge1}H^n(A,A)$ is zero. 
\end{proof} 
 
\ 
 
\subsection{A distinguished example} 
\label{subsec:example} 
 
We exhibit a non cycle TQA with non trivial cohomology ring. 
It turns out that this example is a fundamental piece to understand 
and describe the full structure of the cohomology ring of any TQA. 
This will be done in a forthcomming paper. 
 
Let $\Delta$ be the following quiver 
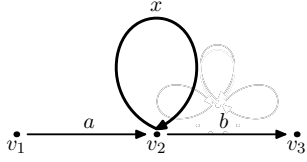
\begin{center} 
   \includegraphics[width=40mm]{quiver.2} 
\end{center} 
and let $A={\mathrm{k}}\Delta/(\Delta_N)$ be a TQA with $N\ge3$. 
It is clear that $\Delta_0\parallel\Delta_0=\{(v_1,v_1),(v_2,v_2),(v_3,v_3)\}$ and 
it is not difficult to see  that 
for $j>0$ or $M>0$ 
\[ 
\Delta_j\parallel\Delta_M= 
\left\{\; 
(ax^{j-1}, ax^{M-1})\;,\; 
(x^{j}, x^{M})\;,\; 
(x^{j-1}b, x^{M-1}b)\;,\; 
(ax^{j-2}b, ax^{M-2}b) 
\;\right\} 
\] 
with the convention that $x^0=v_2$ and a pair containing $x^m$ with $m<0$ do not appear. 
Observe that for this quiver we have that for any $(\alpha,\pi)\in\Delta_j\parallel\Delta_M$ 
the path $\pi$ is determined by $\alpha$. 
Thus, in order to have a clearer notation we shall denote the pair $(\alpha,\pi)$ 
by $\alpha$. 
Now, using this notation, we have 
\begin{align*} 
\Delta_0\parallel\Delta_0&= 
\left\{\; 
v_1,\; 
v_2,\; 
v_3\;\right\},\text{ and} \\ 
\Delta_j\parallel\Delta_M&= 
\left\{\; 
ax^{j-1},\; 
x^{j},\; 
x^{j-1}b,\; 
ax^{j-2}b 
\;\right\},\text{ for $j>0$ or $M>0$}. 
\end{align*} 
The associated matrices of the differentials of the complex of Theorem \ref{Thm:Lo} 
are described below: 
 
\noindent 
for $k=0$ 
\begin{center} 
\setlength{\unitlength}{0.94cm} 
\begin{picture}(10,7)(1.2,-7.4) 
  \put(-.3,-1){$\Hom_{\Delta_0}({\mathrm{k}}\Delta_{0},A)$} 
  \put(2.6,-.9){\vector(1,0){1.3}} 
  \put(3.1,-0.75){$\scriptstyle d_{0}$} 
  \put(4.1,-1){$\Hom_{\Delta_0}( {\mathrm{k}}\Delta_{1},A)$} 
  \put(7,-.9){\vector(1,0){1.8}} 
  \put(7.8,-0.75){$\scriptstyle d_{1}$} 
  \put(9.0,-1){$\Hom_{\Delta_0}( {\mathrm{k}}\Delta_{2},A)$} 
\put(-.5,-1.3){\line(1,0){12.4}} 
  \put(0,-2){${\mathrm{k}}\;\Delta_0\parallel\Delta_{0}$} 
  \put(4.8,-2){${\mathrm{k}}\;\Delta_0\parallel\Delta_{1}$} 
  \put(9.5,-2){${\mathrm{k}}\;\Delta_0\parallel\Delta_{2}$} 
  \put(0,-3){${\mathrm{k}}\;\Delta_1\parallel\Delta_{0}$} 
  \put(4.8,-3){${\mathrm{k}}\;\Delta_1\parallel\Delta_{1}$} 
  \put(9.5,-3){${\mathrm{k}}\;\Delta_1\parallel\Delta_{2}$} 
  \put(0,-4){${\mathrm{k}}\;\Delta_2\parallel\Delta_{0}$} 
  \put(4.8,-4){${\mathrm{k}}\;\Delta_2\parallel\Delta_{1}$} 
  \put(9.5,-4){${\mathrm{k}}\;\Delta_2\parallel\Delta_{2}$} 
  \put(1,-4.8){$\cdot$} 
  \put(5.7,-4.8){$\cdot$} 
  \put(10.5,-4.8){$\cdot$} 
  \put(1,-5){$\cdot$} 
  \put(5.7,-5){$\cdot$} 
  \put(10.5,-5){$\cdot$} 
  \put(1,-5.2){$\cdot$} 
  \put(5.7,-5.2){$\cdot$} 
  \put(10.5,-5.2){$\cdot$} 
  \put(0,-6){${\mathrm{k}}\;\Delta_{N-2}\parallel\Delta_{0}$} 
  \put(4.6,-6){${\mathrm{k}}\;\Delta_{N-2}\parallel\Delta_{1}$} 
  \put(9.5,-6){${\mathrm{k}}\;\Delta_{N-2}\parallel\Delta_{2}$} 
  \put(0,-7){${\mathrm{k}}\;\Delta_{N-1}\parallel\Delta_{0}$} 
  \put(4.6,-7){${\mathrm{k}}\;\Delta_{N-1}\parallel\Delta_{1}$} 
  \put(9.5,-7){${\mathrm{k}}\;\Delta_{N-1}\parallel\Delta_{2}$} 
  \put(2.4,-2.2){\vector(3,-1){1.8}}\put(3.2,-2.3){$\scriptstyle D_0^{0}$} 
  \put(2.4,-3.2){\vector(3,-1){1.8}}\put(3.2,-3.3){$\scriptstyle D_1^{0}$} 
  \put(3.3,-4.8){$\cdot$} 
  \put(3.3,-5.0){$\cdot$} 
  \put(3.3,-5.2){$\cdot$} 
  \put(2.4,-6.2){\vector(3,-1){1.8}}\put(3.2,-6.3){$\scriptstyle D_{N-2}^{0}$} 
  \put(8.3,-5){$\scriptstyle D_{0}^{1}$} 
  \put(7.1,-1.9){\line(1,0){.7}} 
  \put(7.8,-2.15){\oval(0.5,0.5)[tr]} 
  \put(8.05,-2.15){\line(0,-1){4.5}} 
  \put(8.3,-6.65){\oval(0.5,0.5)[bl]} 
  \put(8.3,-6.9){\vector(1,0){.7}} 
%
\end{picture} 
\end{center} 
we have 
\[ 
[D_0^{0}]= 
\scriptsize{ 
\left[\!\begin{array}{rrr} 
-1& 1& 0   \\ 
 0& 0& 0   \\ 
 0&-1& 1 
\end{array}\!\right]}\!;\quad 
[D_j^{0}]= 
\scriptsize{ 
\left[\begin{array}{r} 
 1   \\ 
 0   \\ 
-1 
\end{array}\right]}\!,\normalsize\text{ $2\le j\le N\!-\!2$} 
;\quad 
[D_0^{1}]= 
\scriptsize{ 
\left[\!\begin{array}{c} 
N-1   \\ 
N     \\ 
N-1   \\ 
N-2 
\end{array}\!\right]}\!, 
\] 
and for $k\ge1$ 
\begin{center} 
\setlength{\unitlength}{0.94cm} 
\begin{picture}(10,7)(1.2,-7.4) 
  \put(-.4,-1){$\Hom_{\Delta_0^e}({\mathrm{k}}\Delta_{kN},A)$} 
  \put(2.8,-.9){\vector(1,0){1.1}} 
  \put(3.1,-0.75){$\scriptstyle d_{2k}$} 
  \put(4.2,-1){$\Hom_{\Delta_0^e}( {\mathrm{k}}\Delta_{kN+1},A)$} 
  \put(7.8,-.9){\vector(1,0){1.1}} 
  \put(8,-0.75){$\scriptstyle d_{2k+1}$} 
  \put(9.1,-1){$\Hom_{\Delta_0^e}( {\mathrm{k}}\Delta_{(k+1)N},A)$} 
\put(-.5,-1.3){\line(1,0){13.3}} 
  \put(0,-2){${\mathrm{k}}\;\Delta_0\parallel\Delta_{ kN}$} 
  \put(4.6,-2){${\mathrm{k}}\;\Delta_0\parallel\Delta_{ kN+1}$} 
  \put(9.5,-2){${\mathrm{k}}\;\Delta_0\parallel\Delta_{ (k+1)N}$} 
  \put(0,-3){${\mathrm{k}}\;\Delta_1\parallel\Delta_{ kN}$} 
  \put(4.6,-3){${\mathrm{k}}\;\Delta_1\parallel\Delta_{ kN+1}$} 
  \put(9.5,-3){${\mathrm{k}}\;\Delta_1\parallel\Delta_{ (k+1)N}$} 
  \put(0,-4){${\mathrm{k}}\;\Delta_2\parallel\Delta_{ kN}$} 
  \put(4.6,-4){${\mathrm{k}}\;\Delta_2\parallel\Delta_{ kN+1}$} 
  \put(9.5,-4){${\mathrm{k}}\;\Delta_2\parallel\Delta_{ (k+1)N}$} 
  \put(1,-4.8){$\cdot$} 
  \put(6,-4.8){$\cdot$} 
  \put(10.5,-4.8){$\cdot$} 
  \put(1,-5){$\cdot$} 
  \put(6,-5){$\cdot$} 
  \put(10.5,-5){$\cdot$} 
  \put(1,-5.2){$\cdot$} 
  \put(6,-5.2){$\cdot$} 
  \put(10.5,-5.2){$\cdot$} 
  \put(0,-6){${\mathrm{k}}\;\Delta_{N-2}\parallel\Delta_{ kN}$} 
  \put(4.6,-6){${\mathrm{k}}\;\Delta_{N-2}\parallel\Delta_{ kN+1}$} 
  \put(9.5,-6){${\mathrm{k}}\;\Delta_{N-2}\parallel\Delta_{ (k+1)N}$} 
  \put(0,-7){${\mathrm{k}}\;\Delta_{N-1}\parallel\Delta_{ kN}$} 
  \put(4.6,-7){${\mathrm{k}}\;\Delta_{N-1}\parallel\Delta_{ kN+1}$} 
  \put(9.5,-7){${\mathrm{k}}\;\Delta_{N-1}\parallel\Delta_{ (k+1)N}$} 
  \put(2.4,-2.2){\vector(3,-1){1.8}}\put(3.2,-2.3){$\scriptstyle D_0^{2k}$} 
  \put(2.4,-3.2){\vector(3,-1){1.8}}\put(3.2,-3.3){$\scriptstyle D_1^{2k}$} 
  \put(3.3,-4.8){$\cdot$} 
  \put(3.3,-5.0){$\cdot$} 
  \put(3.3,-5.2){$\cdot$} 
  \put(2.4,-6.2){\vector(3,-1){1.8}}\put(3.2,-6.3){$\scriptstyle D_{N-2}^{2k}$} 
  \put(8.5,-5){$\scriptstyle D_{0}^{2k+1}$} 
  \put(7.6,-1.9){\line(1,0){.5}} 
  \put(8.1,-2.15){\oval(0.5,0.5)[tr]} 
  \put(8.35,-2.15){\line(0,-1){4.5}} 
  \put(8.6,-6.65){\oval(0.5,0.5)[bl]} 
  \put(8.6,-6.9){\vector(1,0){.5}} 
%
\end{picture} 
\end{center} 
we have 
\[ 
[D_0^{2k}]= 
\scriptsize{ 
\left[\!\begin{array}{r} 
 1   \\ 
 0   \\ 
-1 
\end{array}\!\right]},\quad 
[D_0^{2k+1}]= 
\scriptsize{ 
\left[\!\begin{array}{c} 
N-1   \\ 
N     \\ 
N-1   \\ 
N-2 
\end{array}\!\right]}, 
\] 
\[ 
[D_1^{2k}]= 
\scriptsize{ 
\left[\begin{array}{rrr} 
-1& 1& 0   \\ 
 0& 0& 0   \\ 
 0&-1& 1   \\ 
-1& 0& 1 
\end{array}\right]},\quad 
[D_j^{2k}]= 
\scriptsize{ 
\left[\begin{array}{rrrr} 
-1& 1& 0& 0   \\ 
 0& 0& 0& 0   \\ 
 0&-1& 1& 0   \\ 
-1& 0& 1& 0 
\end{array}\right]},\normalsize\text{  $2\le j\le N\!-\!2$.} 
\] 
Therefore a basis of the cohomology is described by the following table. 
 
\ 
 
{\noindent\scriptsize 
\begin{tabular}{|c|cccc|} 
  \hline 
         & $\kk\Delta_0\parallel\Delta_M$     & $\kk\Delta_1\parallel\Delta_M$ 
         & $\kk\Delta_2\parallel\Delta_M,\dots,\kk\Delta_{N\!-\!2}$ 
         & $\kk\Delta_{N\!-\!1}\parallel\Delta_M$ 
                                                \rule[-4pt]{0pt}{16pt}\\ 
  \hline 
$H^0(A,A)$    & 1 
            & $\emptyset$     & $\emptyset$     & $x^{n-1}$ 
                                                \rule[-0pt]{0pt}{12pt}   \\ 
$\dim=2$ &&&&                                      \rule[-4pt]{0pt}{16pt}   \\ 
  \hline 
$H^{1}(A,A)$    & $\emptyset$ 
            & 
            $\begin{matrix} 
            x+a \\[1mm] 
            \\ 
            \end{matrix}$ 
            & 
            $\begin{matrix} 
            x^{j}+ax^{j\!-\!1} \\[1mm] 
            ax^{j\!-\!1} 
            \end{matrix}$ 
            & 
            $\begin{matrix} 
            x^{N\!-\!1}+ax^{N\!-\!2} \\[1mm] 
            ax^{N\!-\!2} 
            \end{matrix}$ 
                                                \rule[-6pt]{0pt}{28pt}   \\ 
\tiny $\dim=2N\!-\!3$ &&&&                                      \rule[-10pt]{0pt}{12pt}   \\ 
\tiny coboundaries    & 
            & $a$, $b$ 
            & $ax^{j\!-\!1}\!-\!x^{j\!-\!1}b$ 
            & $ax^{N\!-\!2}\!-\!x^{N\!-\!2}b$ 
                                                \rule[-6pt]{0pt}{16pt} \\ 
  \hline 
$\begin{matrix} 
H^{2k}(A,A) \\ 
k\ge1 
\end{matrix}$ 
    & $\emptyset$ 
            & 
            $\begin{matrix} 
            a+x+b \\[1mm] 
            \\[1mm] 
            \\ 
            \end{matrix}$ 
           & 
            $\begin{matrix} 
            ax^{j\!-\!1}\!+\!x^{j}\!+\!x^{j\!-\!1}b\!+\!ax^{j\!-\!2}b \\[1mm] 
            ax^{j\!-\!2}b\\[1mm] 
            \\ 
            \end{matrix}$ 
            & 
            $\begin{matrix} 
            ax^{N\!-\!2}\!+\!x^{N}\!+\!x^{N\!-\!2}b\!+\!ax^{N\!-\!3}b \\[1mm] 
            ax^{N\!-\!3}b \\[1mm] 
            x^{N\!-\!1} 
            \end{matrix}$ 
                                                \rule[-6pt]{0pt}{32pt}   \\ 
\tiny $\dim=2N\!-\!2$ &&&&                                      \rule[-10pt]{0pt}{12pt}   \\ 
\tiny coboundaries    & 
            & $\emptyset$ 
            & $\emptyset$ 
            & $\begin{matrix} 
              (N\!-\!1)ax^{N\!-\!2}\!+\!Nx^{N\!-\!1}\!+\!\\ 
              \!+\!(N\!-\!1)x^{N\!-\!2}b\!+\!(N\!-\!2)ax^{N\!-\!3}b 
            \end{matrix}$ 
                                                \rule[-6pt]{0pt}{16pt} \\[3mm] 
  \hline 
$\begin{matrix} 
H^{2k+1}(A,A) \\ 
k\ge1 
\end{matrix}$ 
            & $\emptyset$ 
            & 
             $\begin{matrix} 
            x+a \\[1mm] 
            a 
            \end{matrix}$ 
           & 
            $\begin{matrix} 
            x^{j}+ax^{j\!-\!1} \\[1mm] 
            ax^{j\!-\!1} 
            \end{matrix}$ 
            & 
            $\begin{matrix} 
            x^{N\!-\!1}+ax^{N\!-\!2} \\[1mm] 
            ax^{N\!-\!2} 
            \end{matrix}$ 
                                                \rule[-12pt]{0pt}{32pt}   \\ 
\tiny $\dim=2N\!-\!2$ &&&&                                      \rule[-8pt]{0pt}{20pt}   \\ 
\tiny coboundaries   & 
            & $a\!-\!b$ 
            & $\begin{matrix} 
            ax^{j\!-\!1}\!+\!ax^{j\!-\!2}b, \\[1mm] 
            x^{j\!-\!1}b\!+\!ax^{j\!-\!2}b 
            \end{matrix}$ 
            & $\begin{matrix} 
            ax^{N\!-\!2}\!+\!ax^{N\!-\!3}b, \\[1mm] 
            x^{N\!-\!2}b\!+\!ax^{N\!-\!3}b 
            \end{matrix}$ 
                                                \rule[-6pt]{0pt}{16pt} \\[3mm] 
  \hline 
\end{tabular} 
} 
 
\ 
 
The elements of this basis have been chosen so that the cup product becomes more transparent. 
For $n\ge1$ and $1\le j\le N-1$ let $\omega_{n,j}$ be the basis element 
of $\kk\Delta_j\parallel\Delta_M\subset H^{n}(A,A)$ placed at the top of each row in the above table, 
that is 
\[ 
\omega_{n,j}= 
\begin{cases} 
x^{j}+ax^{j\!-\!1},                              &\text{ if $n$ is odd}; \\ 
ax^{j\!-\!1}+x^{j}+x^{j\!-\!1}b+ax^{j\!-\!2}b,   &\text{ if $n$ is even}. 
\end{cases} 
\] 
Note that $\omega_{2k,j}$ is a sum of two different medal cohomology classes 
(see Definition \ref{Def:medal class}) $\bar M_1+\bar M_2$ where 
$M_1=\{ax^{j\!-\!1},x^{j},x^{j\!-\!1}b\}$ and 
$M_2=\{ax^{j\!-\!2}b\}$. 
From Theorem \ref{Thm:YonedaCohomology} it follows that 
\[ 
\omega_{n_1,j_1}\cup\omega_{n_2,j_2}= 
\begin{cases} 
\omega_{n_1+n_2,j_1+j_2}, &\text{ if $n_1$ or $n_2$ is even and $j_1+j_2<N$}; \\ 
0,&\text{  otherwise}. 
\end{cases} 
\] 
 
\ 
 
\subsection{Non zero cohomology classes in the bar complex} 
\label{subsec:Constructing} 
 
In this subsection we use the comparison morphism to construct explicit 
non zero cohomology classes in the bar complex. 
In the first example we consider the group $H^{2k}_{N-1}$ 
of any $N$ TQA and in the last one we give a full description of the 
cohomology ring of truncated polynomial algebras in one variable.

\subsubsection{Non zero cohomology classes in $H^{2k}_{N-1}$}\ 
 
Recall that $H^{2k}_{N-1}$ is the cokernel of the injective map 
\begin{align*} 
D_0^{2k-1}:\Delta_0\parallel\Delta_{(k-1)N+1}&\to\Delta_{N-1}\parallel\Delta_{ kN} \\ 
(v,\pi)&\mapsto\sum_{ab\in\Delta_{N-1}}(avb,a\pi b), 
\end{align*} 
(see Theorem \ref{Thm:Lo} and Remark \ref{rmk:bigraded}). 
In particular, a pair of parallel paths 
$(\beta,\tau)\in \Delta_{N-1}\parallel\Delta_{ kN}$ with the 
property that it neither start together nor end together is not in the image of 
$D_0^{2k-1}$ and hence corresponds to a non zero cohomology class. 
 
Assume that there exists such a pair $(\beta,\tau)\in \Delta_{N-1}\parallel\Delta_{ kN}$. 
According to the identification \eqref{HomP1} it corresponds to the 
element $g_{(\beta,\tau)}\in\Hom_{\Delta_{0}^e}({\mathrm{k}}\Delta_{kN},A)$ given by 
\[ 
g_{(\beta,\tau)}(\pi)= 
\begin{cases} 
  \beta,&\text{if $\pi=\tau$}; \\ 
  0&\text{otherwise}. 
\end{cases} 
\] 
It is straightforward to see that 
\[ 
f_{(\beta,\tau)}=g_{(\beta,\tau)}\circ{\bf G}\in\Hom_{A^e}(A\otimes A^{\otimes 2k}\otimes A,A) 
\simeq\Hom_{\mathrm{k}}(A^{\otimes 2k},A) 
\] 
is given by 
\begin{align*} 
f_{(\beta,\tau)}(1[\alpha_1|\dots|\alpha_{2k}]1)= 
\begin{cases} 
\beta,& 
\text{if }\alpha_{2i-1}\alpha_{2i}=0\text{ in $A$ for }i=1\dots k \\[-1mm] 
 &\text{and $\alpha_1\dots\alpha_{2k}=\tau$ in $\kk\Delta$};\\[2mm] 
0&\text{otherwise}. 
\end{cases} 
\end{align*}

\subsubsection{Truncated polynomial algebra in one variable}\ 
 
If $A=\kk[x]/(x^N)$ we have 
\[ 
{\bf P}^*_{2k,i}={\mathrm{k}}(x^i,x^{kN}) \qquad\text{and}\qquad 
{\bf P}^*_{2k+1,i}={\mathrm{k}}(x^i,x^{kN+1}) 
\] 
and the only non zero differentials are $D_0^{2k+1}$ for all $k$. 
Thus 
\begin{align*} 
H^{2k}(A,A)&={\mathrm{k}}\{(x^0,x^{kN}),(x^1,x^{kN}),\dots,(x^{N-2},x^{kN})\} \\ 
\qquad\text{and}\qquad 
H^{2k+1}(A,A)&={\mathrm{k}}\{(x^1,x^{kN+1}),(x^2,x^{kN+1}),\dots,(x^{N-1},x^{kN+1})\}. 
\end{align*} 
Using the comparison morphism ${\bf G}$ it is not difficult to give a basis of 
the cohomology in the (reduced) bar resolution. 
Indeed, let 
\[ 
q=1[x^{r_1}|\dots|x^{r_{n}}]1\in {\bf Q}_n= 
A\otimes_{\Delta_0}A_+^{\otimes^n_{\Delta_0}}\otimes_{\Delta_0}A,\qquad 
r_i>0; 
\] 
and let 
\begin{align*} 
f_{2k,i}(1[x^{r_1}|\dots|x^{r_{2k}}]1)= 
\begin{cases} 
x^{i+\sum r_j-kN},& 
\text{if }r_{2j-1}+r_{2j}\ge N\text{ for }j=1\dots k; \\[2mm] 
0, &\text{otherwise}; 
\end{cases} 
\end{align*} 
\begin{align*} 
f_{2k+1,i}(1[x^{r_1}|\dots|x^{r_{2k+1}}]1)= 
\begin{cases} 
x^{i+\sum r_j-kN-1},& 
\text{if }r_{2j}+r_{2j+1}\ge N\text{ for }j=1\dots k; \\[2mm] 
0, &\text{otherwise}. 
\end{cases} 
\end{align*} 
Then 
\begin{align*} 
H^{2k}(A,A)&={\mathrm{k}}\{f_{2k,0},f_{2k,1},\dots,f_{2k,N-2}\} \\ 
\qquad\text{and}\qquad 
H^{2k+1}(A,A)&={\mathrm{k}}\{f_{2k+1,1},f_{2k+1,2},\dots,f_{2k+1,N-1}\}. 
\end{align*} 
The cup product is given by 
\[ 
f_{m,i}\cup f_{n,j}= 
\begin{cases} 
f_{m+n,i+j} & \text{if either $m$ or $n$ is even and $i+j<N$};\\ 
0, &\text{otherwise}; 
\end{cases} 
\] 
and $\{f_{0,1},f_{1,1},f_{2,0}\}$ is a set of generators of $H^{*}(A,A)$. 
 

\end{document}